\documentclass[12pt]{amsart}
\usepackage[]{amscd,amsfonts,amsmath,amsxtra,amssymb}
\usepackage{hyperref}
\usepackage{graphics}
\usepackage[all]{xy}
\usepackage[left]{lineno}
\newtheorem{sub}{}[section]
\newtheorem{subsub}{}[sub]

\def\HHom{\mathop{\mathcal Hom}\nolimits}

\def\EEnd{\mathop{\mathcal End}\nolimits}

\def\QQuot{\mathop{\mathcal Quot}\nolimits}

\def\hfl#1#2{\smash{\mathop{\ \hbox to 12mm{\rightarrowfill}}
\limits^{\scriptstyle#1}_{\scriptstyle#2} \ }}
\def\hflb#1#2{\smash{\mathop{\hbox to 12mm{\leftarrowfill}}
\limits^{\scriptstyle#1}_{\scriptstyle#2}}}

\def\og{\leavevmode\raise.3ex\hbox{$\scriptscriptstyle\langle\!\langle$}}
\def\fg{\leavevmode\raise.3ex\hbox{$\scriptscriptstyle\,\rangle\!\rangle$}}

\def\Ssect#1#2{\pagebreak[3]\begin{sub}\label{#2}{\sc\small\small
#1}\rm\medskip}

\def\xmat#1{\[\xymatrix{#1}\]}
\def\flinc{\ar@{^{(}->}}
\def\fleq{\ar@{=}}
\def\flon{\ar@{->>}}
\def\fmaps{\ar@{|-{>}}}

\def\SS{\mathop{\mathcal S}\nolimits}

\newcommand{\N}{{\mathbb N}}

\newcommand{\Z}{{\mathbb Z}}

\newcommand{\C}{{\mathbb C}}

\renewcommand{\P}{{\mathbb P}}

\newcommand{\kh}{{\mathcal H}}

\newcommand{\km}{{\mathcal M}}
\newcommand{\kn}{{\mathcal N}}
\newcommand{\ko}{{\mathcal O}}

\newcommand{\kq}{{\mathcal Q}}

\newcommand\bigzero{\makebox(0,0){\text{\huge0}}}

\begin{document}
\def\refname{R\'ef\'erences}
\def\contentsname{ }
\def\proofname{D\'emonstration}
\def\abstractname{R\'esum\'e}

\author{Mohamed Bahtiti}
\address{Institut de Math\'{e}matiques de Jussieu,
 Case 247, 4 place Jussieu,\newline
F-75252 Paris, France}
\email{mohamed.bahtiti@imj-prg.fr}
\title{ Fibr\'e vectoriel de 0-corr\'elation pond\'er\'e sur~l'espace~$\P^{2n+1}$ }
 
\maketitle

\begin{abstract}
Nous \'etudions dans cet article une nouvelle famille de fibr\'es vectoriels symplectiques alg\'ebriques stables de rang $2n$ sur l'espace projectif complexe $\P^{2n+1}$ dont le fibr\'e de corr\'elation nulle classique fait partie. Nous montrons que cette famille est invariante par rapport aux d\'eformations miniversales. Nous \'etudions \'egalement les conditions cohomologiques suffisantes pour qu'un fibr\'e vectoriel symplectique sur une vari\'et\'e projective soit stable. 
$$$$

{\scriptsize ABSTRACT}. We study in this paper a new family of stable algebraic symplectic vector bundles of rank $ 2n $ on the complex projective space $\P^{2n+1}$ whose classical null correlation bundles belongs. We show that these bundles are invariant under a miniversal deformation. We also study the sufficient cohomological conditions for a symplectic vector bundle on a projective variety to be stable.
 
\end{abstract}

{\scriptsize {\em Date}: August, 2015.
{\em Mots-cl\'es}. fibr\'e de corr\'elation nulle, stabilit\'e, d\'eformation miniversale, espace de Kuranishi.}

\vspace{1.8cm}

\section{\large\bf Introduction} 
 
\vspace{1cm}

Les fibr\'es vectoriels alg\'ebriques non-d\'ecomposables connus de rang $n-1$ sur l'espace projectif complexe $\P^{n}$ pour $n \geq 6$ sont rares. Les familles de fibr\'es vectoriels seulement connues sont la famille de fibr\'es de Tango pond\'er\'es de rang $n-1$ \cite{ca} et celle de fibr\'es instantons de rang $n-1$ pour $n$ impair \cite{ok-sp}. 

La famille de fibr\'es instantons pr\'esente un lien important entre la g\'eom\'etrie alg\'ebrique et la physique math\'ematique, en particulier la th\'eorie de Yang-Mills. Cette famille de fibr\'es a \'et\'e construite sur $\P^{3}$ par Atiyah, Drenfeld, Hitchin et Manin \cite{atdrhima} et correspond, gr\^ace \`a la \mbox{correspondance} de Penrose-Ward, aux solutions auto-duales (instantons) des  $SU(2)$-\'equations de Yang-Mills sur la sph\`ere euclidienne $S^{4}$ \cite{atdrhima,dove, ati, atwa}. Suite \`a la g\'en\'eralisation de la correspondance de Penrose-Ward sur l'espace projectif complexe de dimension impaire par Salamon \cite{sa}, la famille de fibr\'es instantons a \'et\'e g\'en\'eralis\'ee sur $\P^{2n+1}$ par Okonek et Spindler \cite{ok-sp}. Ensuite Spindler et Trautmann ont \'etudi\'e la famille de fibr\'es instantons sp\'eciaux \cite{sp-tr}.

Les fibr\'es instantons sp\'eciaux de nombre quantique $1$ sont des fibr\'es de corr\'elation nulle \mbox{classiques} sur $\P^{2n+1}$ et sont la cohomologie de monades de type
$$0\longrightarrow \ko_{\P^{2n+1}}(-1)\stackrel{B}{\longrightarrow} \C^{2n+2}\otimes \ko_{\P^{2n+1}} \stackrel{A}{\longrightarrow}\ko_{\P^{2n+1}}(1) \longrightarrow 0,$$

o\`u 
$$A =
\left(
\begin{array}{ccccccc}
y_{n} & \ldots & y_{0} & ; & -x_{n}& \ldots &-x_{0} \\
\end{array}
\right), \hspace{0.2 cm}  B =   ^{T}\left(
\begin{array}{ccccccc}
x_{0} & \ldots & x_{n} &;&y_{0}&  \ldots  &y_{n} \\
\end{array}
\right), $$

$x_{0},x_{1}, \ldots , x_{n}, y_{0},y_{1}, \ldots  ,y_{n}$ \'etant des formes lin\'eaires sans z\'ero commun sur $\P^{2n+1}$.

En particulier nous nous int\'eressons \`a la g\'en\'eralisation de ces fibr\'es qui ont \'et\'e \'etudi\'es sur $\P^3$ par Ein \cite{ei}, sur $\P^{5}$ par Ancona et Ottaviani pour la premi\`ere classe de Chern $c_{1}=0$  \cite{an-ot93} et qui ont \'et\'e d\'efinis sur $\P^{2n+1}$ par Migliore, Nagel et Peterson \cite{mi-na-pe}. Plus pr\'ecis\'ement, soient $\gamma >0$ et $\zeta=0,1$ et $\lambda_{i}$ des entiers naturels pour $i=0,1,2,  \ldots  ,n$ tels que 
$$\gamma -\zeta>\lambda_{n}\geq \lambda_{n-1}\geq  \ldots   \geq\lambda_{0}\geq 0.$$

Soient $g_{i}\in H^{0}(\P^{2n+1},\ko_{\P^{2n+1}}(\gamma -\lambda_{i}-\zeta))$ et $f_{i}\in H^{0}(\P^{2n+1},\ko_{\P^{2n+1}}(\gamma +\lambda_{n-i}))$ des~formes homog\`enes sur $\P^{2n+1}$ sans z\'ero commun sur $\P^{2n+1}$. Les fibr\'es de corr\'elation nulle pond\'er\'es (les~ fibr\'es de $0$-corr\'elation pond\'er\'es) sont la cohomologie de monades de type
$$0\longrightarrow \ko_{\P^{2n+1}}(-\gamma -\zeta)\stackrel{B}{\longrightarrow} \kh_{\zeta} \stackrel{A}{\longrightarrow}\ko_{\P^{2n+1}}(\gamma) \longrightarrow 0,$$

o\`u
$$A =
\left(
\begin{array}{ccccccc}
g_{n} & \ldots & g_{0} & ; & -f_{n}& \ldots &-f_{0} \\
\end{array}
\right), \hspace{0.2 cm}  B =   ^{T}\left(
\begin{array}{ccccccc}
f_{0} & \ldots & f_{n} &;&g_{0}&  \ldots  &g_{n} \\
\end{array}
\right),$$

et
$$\kh_{\zeta}:= \bigoplus_{i=0}^{n}\left(    \ko_{\P^{2n+1}}(\lambda_{n-i})\right) \oplus \bigoplus_{i=0}^{n}\left( \ko_{\P^{2n+1}}(-\lambda_{i}-\zeta)  \right).  $$

La cohomologie $\kn_{\zeta}$ de la monade pr\'ec\'edente est un fibr\'e vectoriel symplectique normalis\'e de rang
$2n$ sur $\P^{2n+1}$ pour laquelle les conditions suivantes sont \'equivalentes (th\'eor\`eme \ref{3.7})
 
I-  $\gamma -\zeta n > \sum_{i=0}^{n}\lambda_{i}$.

II- $\kn_{\zeta}$ est stable.

III- $\kn_{\zeta}$ est simple.

Les d\'eformations miniversales d'un tel fibr\'e $\kn_{\zeta}$ sont encore des fibr\'es de 0-corr\'elation pond\'er\'es sur $\P^{2n+1}$ et l'espace de Kuranishi du fibr\'e $\kn_{\zeta}$ est lisse au point correspondant de $\kn_{\zeta}$ (th\'eor\`eme \ref{4.8}). Les conditions cohomologiques suffisantes pour qu'un fibr\'e vectoriel symplectique sur une vari\'et\'e projective soit stable sont \'enonc\'ees dans le th\'eor\`eme\ref{0.2.5}.

Je tiens \`a exprimer ma gratitude au directeur de ma th\`ese M. J.-M. Dr\'ezet et au professeur \mbox{G. Ottaviani} pour des discussions utiles. Je remercie \'egalement toutes les personnes qui ont contribu\'e \`a m'aider \`a r\'ealiser mes travaux. Cet article fait partie de ma th\`ese.

\newpage

\section{ \large\bf Pr\'eliminaires} 

\vspace{1cm}

\subsection{ D\'efinition}\label{0.2.1} 

On d\'efinit {\em la r\'esolution de Koszul g\'en\'eralis\'ee} de la suite exacte des fibr\'es
$$0\longrightarrow A \longrightarrow B \longrightarrow F \longrightarrow 0$$

par
$$0\longrightarrow \SS^{i} A \longrightarrow \SS^{i-1} A \otimes  B 
\longrightarrow \SS^{i-2} A\otimes \bigwedge^{2}  B \longrightarrow \SS^{i-3} 
A\otimes \bigwedge^{3}  B \longrightarrow \SS^{i-4} A\otimes \bigwedge^{4}  B 
\longrightarrow \ldots$$
$$\hspace{4 cm} \ldots \longrightarrow \SS^{2} A\otimes \bigwedge^{i-2}  B 
\longrightarrow  A\otimes \bigwedge^{i-1}  B \longrightarrow  \bigwedge^{i}  B 
\longrightarrow \bigwedge^{i} F \longrightarrow 0$$

Soient $rg(F)=r$, et $c_{1}(F)=c$. On a $\bigwedge^{r-i} F^{*}=\bigwedge^{r} F^{*}\otimes \bigwedge^{i} F=\bigwedge^{i} F(-c)$. En tensorisant la r\'esolution de Koszul g\'en\'eralis\'ee par $\bigwedge^{r}F^{*}$, on obtient une r\'esolution pour $\bigwedge^{r-i} F^{*}$
$$0\longrightarrow \SS^{i} A (-c)  \longrightarrow 
\SS^{i-1} A \otimes  B(-c) \longrightarrow \SS^{i-2} 
A\otimes \bigwedge^{2}  B(-c) \longrightarrow$$
$$ \SS^{i-3} A\otimes \bigwedge^{3}  B(-c) 
\longrightarrow \SS^{i-4} A\otimes \bigwedge^{4}  B(-c) 
\longrightarrow \ldots$$
$$\ldots \longrightarrow \SS^{2} A\otimes \bigwedge^{i-2}  B (-c) \longrightarrow  A\otimes \bigwedge^{i-1}  B(-c) 
\longrightarrow  \bigwedge^{i}  B  (-c)\longrightarrow 
\bigwedge^{r-i} F^{*} \longrightarrow 0$$

\subsection{\bf D\'efinition}\label{0.2.2} 

{\em  Soit $E$ un fibr\'e vectoriel de rang $2r$ sur une vari\'et\'e projective $X$. On dit que le fibr\'e $E$ est  symplectique si et seulement s'il existe un entier $b \in \Z$ et un isomorphisme de fibr\'es 
$$\varphi :E \longrightarrow E^{*}(b)\hspace{0.2 cm} tel \hspace{0.2 cm} que\hspace{0.2 cm} \varphi^{*}=-\varphi.$$

Autrement dit, le fibr\'e  $E$ est symplectique si et seulement s'il existe une forme symplectique non-d\'eg\'en\'er\'ee $\varpi \in H^{0}(\bigwedge^{2}E (-b))$}. Dans ce cas, la premi\`ere classe de Chern de $E$ est $c_{1}(E)=br$.

Nous allons d\'evelopper les conditions cohomologiques suffisantes pour qu'un fibr\'e vectoriel symplectique de rang $2r$ sur une vari\'et\'e projective $X$ soit stable. Le lemme suivant est une g\'en\'eralisation de (\cite{an-ot94}, Lemme 1.10).

\subsection{\bf Lemme }\label{0.2.3} 
{\em Soit $E$ un fibr\'e vectoriel symplectique de rang $2r$ \ref{0.2.2} sur une vari\'et\'e projective $X$. Pour tout $1\leq j\leq r-1$, le fibr\'e $\ko_{X}$ est un suppl\'ementaire dans le fibr\'e $\bigwedge^{2j}E(-bj)$ et le fibr\'e $E$ est un suppl\'ementaire dans le fibr\'e $\bigwedge^{2j+1}E(-bj)$. Autrement dit, on a
$$\bigwedge^{2j}E(-bj)\simeq \ko_{X}\oplus B_{j}$$

$$\bigwedge^{2j+1}E(-bj)\simeq E \oplus A_{j} .$$}

\begin{proof}
Pour tout $i\leq r$. On a une forme non-d\'eg\'en\'er\'ee $\varpi \in H^{0}(\bigwedge^{2}E (-b))$, donc on peut d\'efinir un morphisme injectif $\Phi$ qui est localement donn\'e par
$$\Phi: \bigwedge^{i-2}E \longrightarrow \bigwedge^{i}E(-b)$$

$$\hspace{1.1 cm}e_{1}\wedge \ldots \wedge e_{i-2} \longmapsto e_{1}\wedge \ldots \wedge e_{i-2}\wedge \varpi. $$

Alors pour $i=2j$, le fibr\'e $\ko_{X}$ est un suppl\'ementaire dans le fibr\'e $\bigwedge^{2j}E(-jb)$.
Pour $i=2j+1$, le fibr\'e $E$ est aussi un suppl\'ementaire dans le fibr\'e $\bigwedge^{2j+1}E(-jb)$.   

\end{proof}

\subsection{\bf Remarque}\label{0.2.4 } 

Pour tout fibr\'e vectoriel, le fibr\'e $ \bigwedge^{j+1}E$ est un suppl\'ementaire dans le fibr\'e $E\otimes \bigwedge^{j}E$, o\`u $j\geq 0$. C'est-\`a-dire, on a un morphisme injectif localement donn\'e par
$$\Phi: \bigwedge^{j+1}E \longrightarrow E\otimes \bigwedge^{j}E\hspace{4 cm}$$

$$\hspace{2.2 cm} e_{1} \wedge \ldots \wedge e_{j+1} \longmapsto \frac{1}{j+1} \sum_{i=1}^{j+1}(-1)^{j-i+1} (e_{1}\wedge \ldots \wedge \widehat{e_{i} } \wedge \ldots \wedge e_{j+1})\otimes e_{i}. $$

Le th\'eor\`eme suivant est une g\'en\'eralisation de (\cite{an-ot94}, th\'eor\`eme 3.5).\\

\subsection{\bf Th\'eor\`eme}\label{0.2.5} 
{\em Soient $E$ un fibr\'e vectoriel symplectique de rang $2r$, $E \simeq E^{*}(b)$ o\`u $b\in \Z$, sur une vari\'et\'e projective $X$ avec $Pic(X)=\Z$, et $1\leq 2j+1 \leq r$ un entier. On a

1- Soient $b\geq 0$, et

 \hspace{0.5 cm} I) $h^{0}(\bigwedge^{2j+1}E)=0$,

 \hspace{0.5 cm} II) $h^{0}(E(-b(j+1))\otimes \bigwedge^{2j+1}E)=1$.

 \hspace{0.5 cm} Alors $E$ est stable (au sens de Mumford-Takemato).

2- Soient $b\leq 0$, et

 \hspace{0.5 cm} I) $h^{0}(\bigwedge^{2j+1}E^{*})=0$,

\hspace{0.5 cm}  II) $h^{0}(E^{*}(b(j+1)) \otimes \bigwedge^{2j+1}E^{*})=1$.

\hspace{0.5 cm}  Alors $E$ est stable (au sens de Mumford-Takemato)}.

\begin{proof}
1- Pour $b\geq 0$. Premi\`erement, supposons que l'on ait  $Z_{1}$ un sous-faisceau de fibr\'e $E$ de rang $2t_{1}+1$ qui d\'estabilise le fibr\'e $E$. D'apr\`es la d\'efinition de la stabilit\'e on a $c_{1}(Z_{1})=m_{1}\geq 0$, alors on a la suite exacte suivante
$$0\longrightarrow Z_{1}\longrightarrow E \longrightarrow G \longrightarrow 0$$

o\`u $G$ est un faisceau sans torsion. En appliquant la r\'esolution de Koszul \ref{0.2.1}, on a 
$$0\longrightarrow \bigwedge^{2t_{1}+1}Z_{1} \longrightarrow \bigwedge^{2t_{1}+1}E ,$$ 

qui donne le morphisme
$$0\longrightarrow \ko_{X}(m_{1}) \longrightarrow \bigwedge^{2t_{1}+1}E, $$

on obtient

\xmat{ 0 \ar[rr] &  & \ko_{X}  \ar[rr] \ar[rd]^{\varphi} &  & \bigwedge^{2t_{1}+1}E(-m_{1})\ar[ld] \\
&  & &    \bigwedge^{2t_{1}+1}E }

Donc le morphisme $\varphi$ est une section de fibr\'e $\bigwedge^{2t_{1}+1}E$, 
ce qui contredit $h^{0}(\bigwedge^{2t_{1}+1}E)=0$.

Deuxi\`emement. Supposons que l'on ait $Z$ un sous-faisceau de fibr\'e $E$ de rang $2t+2$ et $c_{1}(Z)=m$ qui d\'estabilise le fibr\'e $E$, alors on a la suite exacte suivante
$$(*) \hspace{2 cm} 0\longrightarrow Z\longrightarrow E \longrightarrow G \longrightarrow 0$$

o\`u $G$ est un faisceau sans torsion. En appliquant la r\'esolution de Koszul \ref{0.2.1}, on a 
$$0\longrightarrow \bigwedge^{2t+2}Z \stackrel{f}{\longrightarrow} \bigwedge^{2t+2}E \stackrel{g}{\longrightarrow} G \otimes \bigwedge^{2t+1}E \longrightarrow \ldots $$ 

qui donne  
$$0\longrightarrow \ko_{X}(m-b(t+1)) \stackrel{f}{\longrightarrow} \bigwedge^{2t+2}E(-b(t+1)) \stackrel{g}{\longrightarrow} G(-b(t+1)) \otimes \bigwedge^{2t+1}E  \longrightarrow \ldots$$
 
o\`u on a $g.f=0$. En tensorisant la suite exacte (*) par le fibr\'e $\bigwedge^{2t+1}E(-b(t+1)) $, on obtient la suite exacte suivante
$$(**)  \hspace{2 cm} 0\longrightarrow Z\otimes \bigwedge^{2t+1}E (-b(t+1)) \longrightarrow E(-b(t+1))\otimes \bigwedge^{2t+1}E  \stackrel{h}{\longrightarrow} G(-b(t+1))\otimes \bigwedge^{2t+1}E  \longrightarrow 0 .$$
$$$$

Mais d'apr\`es la remarque \ref{0.2.4 } , on a 
$$ \bigwedge^{2t+2}E(-b(t+1)) \subset E(-b(t+1))\otimes \bigwedge^{2t+1}E ,$$

alors on obtient que $h|_{\bigwedge^{2t+2}E(-b(t+1))}=g$. D'apr\`es le lemme \ref{0.2.3} , on obtient   
$$\ko_{X}\oplus B_{t}\cong\bigwedge^{2t+2}E(-b(t+1)) \subset E(-b(t+1))\otimes \bigwedge^{2t+1}E,$$

qui donne, d'apr\`es (II), $h^{0}(B_{t})=0$. Donc on obtient que 
$$f(\ko_{X}(m-b(t+1)))\subset \ko_{X},\hspace{0,2 cm} g(\ko_{X})=0 \hspace{0,2 cm} et \hspace{0,2 cm}  h(\ko_{X})=0.$$

En prenant la cohomologie de la suite exacte (**), on obtient la r\'esolution suivante
$$0\longrightarrow H^{0}( Z(-b)\otimes \bigwedge^{2t+1}E(-bt)) \longrightarrow H^{0}( E(-b)\otimes \bigwedge^{2t+1}E(-bt) )\stackrel{H^{0}(h)}{\longrightarrow} $$
$$H^{0}(G(-b(t+1))\otimes \bigwedge^{2t+1}E) \longrightarrow \ldots $$

Comme 
$$\ko_{X} \subset \bigwedge^{2t+2}E(-b(t+1)) \subset E(-b)\otimes \bigwedge^{2t+1}E(-bt)\simeq \HHom(\bigwedge^{2t+1}E^{*}(bt),E(-b)),$$ 

alors on a un morphisme non nul $ \varphi:\bigwedge^{2t+1}E^{*}(bt) \longrightarrow E (-b) $ correspondant \`a $\ko_{X}$ tel que $ H^{0}(h)(\varphi)=0$. Donc il existe un morphisme surjectif $\psi$ dans $Hom(\bigwedge^{2t+1}E^{*}(bt),Z(-b))$ project\'e sur $\ko_{X}$ tel que

\xmat{  \bigwedge^{2t+1}E^{*}(bt)    \ar[rr]^{\varphi}\ar[rd]^{\psi} &  &  E (-b)\ar[rr] &  & 0  \\
&   Z(-b) \ar[ru]\ar[rd]\\
0\ar[ru] &&0 }

Mais c'est une contradiction au fait que la condition $h^{0}(E(-b(t+1))\otimes \bigwedge^{2t+1}E)=1$ entra\^ine que $\varphi$ est le seul morphisme surjectif non nul project\'e sur la suppl\'ementaire $\ko_{X}$ (voir \ref{0.2.3}). Donc $E$ est stable. 

2- Si $b\leq 0$. On consid\`ere $((E)^{*})^{*}(-b)\simeq E^{*}$. Comme on a
$$h^{0}(\bigwedge^{2j+1}E^{*})=0,\hspace{0.2 cm}  h^{0}(E^{*}(b(j+1)) \otimes \bigwedge^{2j+1}E^{*})=1$$

alors, d'apr\`es (1), $E^{*}$ est stable. Donc $E$ est stable.
 
\end{proof}

\subsection{\bf Remarque}
Soit $E$ un fibr\'e symplectique de rang $2r$ sur une vari\'et\'e projective $X$, \mbox{$E\simeq E^{*}(b)$} pour un $b\in Z $. Alors le fibr\'e normalis\'e de $E$ est un cas parmi les deux cas suivants

- Soient $b=2m-\zeta$ et $\zeta=0,1$. Alors on a $E(-m)\simeq (E(-m))^{*}(-\zeta)$ et le fibr\'e normalis\'e de $E$ est $E(-m)$ avec $c_{1}(E(-m))=-\zeta r$. 

\newpage
 
\section{ \large\bf Fibr\'e vectoriel de 0-corr\'elation pond\'er\'e} 

\vspace{1cm}

Soient $V$ un espace vectoriel complexe de dimension $2n+2$, et $\P^{2n+1}=\P(V)$ l'espace projectif complexe associ\'e dont les points sont les droites de $V$. Soient $\gamma >0$, $n\geq 1$, $\zeta=0,1$ et $\lambda_{i}$ des entiers naturels pour $i=0,1,2,  \ldots  ,n$ tels que 
$$\gamma-\zeta >\lambda_{n}\geq \lambda_{n-1}\geq  \ldots   \geq\lambda_{0}\geq 0.$$

Soient $g_{i}\in H^{0}(\P^{2n+1},\ko_{\P^{2n+1}}(\gamma -\lambda_{i}-\zeta))$ et $f_{i}\in H^{0}(\P^{2n+1},\ko_{\P^{2n+1}}(\gamma +\lambda_{n-i}))$ des~formes homog\`enes sur $\P^{2n+1}$ sans z\'ero commun sur $\P^{2n+1}$. On consid\`ere la $(2n+2)\times 1$-matrice 
$$B =   ^{T}\left(
\begin{array}{ccccccc}
f_{0} & \ldots & f_{n} &;&g_{0}&  \ldots  &g_{n} \\
\end{array}
\right)$$

et la $1\times(2n+2)$-matrice
$$A =
\left(
\begin{array}{ccccccc}
g_{n} & \ldots & g_{0} & ; & -f_{n}& \ldots &-f_{0} \\
\end{array}
\right)$$

qui v\'erifient $A.B=0$. On consid\`ere aussi le fibr\'e 
$$\kh_{\zeta}:= \bigoplus_{i=0}^{n}\left(    \ko_{\P^{2n+1}}(\lambda_{n-i})\right) \oplus \bigoplus_{i=0}^{n}\left( \ko_{\P^{2n+1}}(-\lambda_{i}-\zeta)  \right).  $$

On obtient un isomorphisme canonique $\kh_{\zeta}\simeq \kh_{\zeta}^{*}(-\zeta)$. Tout cela nous donne la monade suivante
$$0\longrightarrow \ko_{\P^{2n+1}}(-\gamma-\zeta)\stackrel{B}{\longrightarrow} \kh_{\zeta} \stackrel{A}{\longrightarrow}\ko_{\P^{2n+1}}(\gamma) \longrightarrow 0.$$

On en d\'eduit les suites exactes suivantes
$$0\longrightarrow \ko_{\P^{2n+1}}(-\gamma-\zeta)\stackrel{B}{\longrightarrow} \kh_{\zeta}\longrightarrow \kq_{\zeta;\gamma;\lambda_{n},\lambda_{n-1},  \ldots  ,\lambda_{0}} \longrightarrow 0,$$

o\`u $\kq_{\zeta;\gamma;\lambda_{n},\lambda_{n-1},  \ldots  ,\lambda_{0}}$ est un fibr\'e vectoriel de rang $2n+1$ sur $\P^{2n+1}$ et
$$0\longrightarrow \kn_{\zeta;\gamma;\lambda_{n},\lambda_{n-1},  \ldots  ,\lambda_{0}}\longrightarrow  \kq_{\zeta;\gamma;\lambda_{n},\lambda_{n-1},  \ldots  ,\lambda_{0}} \stackrel{A}{\longrightarrow} \ko_{\P^{2n+1}}(\gamma) \longrightarrow 0,$$

o\`u $\kn_{\zeta;\gamma;\lambda_{n},\lambda_{n-1},   \ldots  ,\lambda_{0}}$ est un fibr\'e vectoriel de rang $2n$ sur $\P^{2n+1}$. Autrement dit, le fibr\'e vectoriel $\kn_{\zeta;\gamma;\lambda_{n},\lambda_{n-1},   \ldots  ,\lambda_{0}}$ est la cohomologie de la  monade pr\'ec\'edente. Le fibr\'e $\kn_{\zeta;\gamma;\lambda_{n},\lambda_{n-1}, \ldots  ,\lambda_{0}}$ a \'et\'e \'etudi\'e sur $\P^3$ par Ein \cite{ei}, sur $\P^{5}$ par Ancona et Ottaviani pour $\zeta=0$ \cite{an-ot93}, et a \'et\'e d\'efini sur $\P^{2n+1}$ par Migliore, Nagel et Peterson \cite{mi-na-pe}.

\subsection{\bf D\'efinition} \label{2s} 
{\em On appelle le fibr\'e 
$\kq_{\zeta;\gamma;\lambda_{n},\lambda_{n-1},  \ldots  ,\lambda_{0}}$ le fibr\'e de 
quotient pond\'er\'e sur $\P^{2n+1}$ des poids 
$\gamma;\lambda_{n},\lambda_{n-1}, \ldots  ,\lambda_{0}$, et on appelle le 
fibr\'e  $\kn_{\zeta;\gamma;\lambda_{n},\lambda_{n-1}, \ldots  ,\lambda_{0}}$ le 
fibr\'e de 0-corr\'elation pond\'er\'e sur $\P^{2n+1}$ des poids 
$\gamma;\lambda_{n},\lambda_{n-1}, \ldots  ,\lambda_{0}$ avec $c_{1}=-\zeta n$ et d'une charge topologique 
$$ c_{2}=\gamma^{2}- \sum_{i=0}^{n}\lambda_{i}^{2} +\zeta (\gamma +\zeta \dfrac{n.(n-1)}{2} - \sum_{i=0}^{n}\lambda_{i} ).$$}

On va utiliser $\kq_{\zeta}$, $\kn_{\zeta}$ \`a la place de $\kq_{\zeta;\gamma;\lambda_{n},\lambda_{n-1},  \ldots  ,\lambda_{0}}$, $\kn_{\zeta;\gamma;\lambda_{n},\lambda_{n-1},  \ldots  ,\lambda_{0}}$ respectivement. On a donc la monade suivante
\begin{equation}\label{1} 
0\longrightarrow \ko_{\P^{2n+1}}(-\gamma-\zeta)\stackrel{B}{\longrightarrow} \kh_{\zeta} \stackrel{A}{\longrightarrow}\ko_{\P^{2n+1}}(\gamma ) \longrightarrow 0
\end{equation}

et les suites exactes suivantes
\begin{equation} \label{2}
0\longrightarrow \ko_{\P^{2n+1}}(-\gamma-\zeta)\stackrel{B}{\longrightarrow} \kh_{\zeta}\longrightarrow \kq_{\zeta}  \longrightarrow 0,
\end{equation}

et
\begin{equation} \label{3}
0\longrightarrow \kn_{\zeta}\longrightarrow  \kq_{\zeta}  \stackrel{A}{\longrightarrow} \ko_{\P^{2n+1}}(\gamma) \longrightarrow 0.
\end{equation}

On en d\'eduit que $c_{1}(\kq_{\zeta})=\gamma -\zeta n$. Soit $J$ une $(2n+2)\times(2n+2)$-matrice symplectique 
$$J=\left(
\begin{array}{cccccccccc}
  & &  &  &   & & & & & -1 \\
  & &  &\bigzero  &   & & & & .&  \\
  & &  &  &   & & &. & &  \\
  & &  &  &   & &.& & &  \\
  & &  &  &   & -1& & & &  \\
  & &  &  & 1 & & & & &  \\
  & &  &. &   & & & & &  \\
  & &. &  &   & & & & &  \\
  &.&  &  &   & & \bigzero& & &   \\
1 & &  &  &  & & & & &   \\
\end{array}
\right).$$

 Alors on a $J^{2} =-I$ et $^{T}B =-A(-\zeta).J$ et $^{T}A =-J.B(-\zeta)$.
  
\subsection{\bf Proposition }\label{2.3} 
{\em Le fibr\'e de 0-corr\'elation pond\'er\'e $\kn_{\zeta}$ sur $\P^{2n+1}$ est un fibr\'e symplectique}.

\begin{proof}
  
Le fibr\'e $\kn_{\zeta}^{*}$ est la cohomologie de la monade suivante
$$0\longrightarrow \ko_{\P^{2n+1}}(-\gamma )\stackrel{-J.B(-\zeta)}{\longrightarrow} \kh_{\zeta} \stackrel{-A(-\zeta).J}{\longrightarrow}\ko_{\P^{2n+1}}(\gamma+\zeta) \longrightarrow 0$$

o\`u $J^{2} =-I$, $^{T}B =-A(-\zeta).J$ et $^{T}A =-J.B(-\zeta)$. Donc l'isomorphisme de monades suivant

\xmat{ 0 \ar[r] & \ko_{\P^{2n+1}}(-\gamma-\zeta) \ar[rr]^{B } \ar@{=}[d]& &
\kh_{\zeta} \ar[rr]^{A }\ar[d]^{J} & &  \ko_{\P^{2n+1}}(\gamma) \ar[r]\ar[d]^{I(-\zeta)}& 0\\
0 \ar[r] & \ko_{\P^{2n+1}}(-\gamma-\zeta) \ar[rr]^{-J.B(-\zeta) }& &
\kh_{\zeta}^{*}(-\zeta) \ar[rr]^{-A(-\zeta).J }  &&  \ko_{\P^{2n+1}}(\gamma) \ar[r]& 0\\
}

induit un isomorphisme $\varphi :\kn_{\zeta} \longrightarrow \kn_{\zeta}^{*}(-\zeta)$. En transposant le diagramme pr\'ec\'edent, les fl\`eches verticales sont multipli\'ees par $-1$. Donc on obtient $^{T}\varphi (-\zeta)=- \varphi$.

\end{proof}

\subsection{\bf Proposition}\label{2.4} 
{\em Soient $\kq_{\zeta}$ le fibr\'e de quotient pond\'er\'e sur $\P^{2n+1}$ et $\kn_{\zeta}$ 
le fibr\'e de \mbox{$0$-corr\'elation} pond\'er\'e sur $\P^{2n+1}$ qui sont d\'efinis 
par les suites exactes \ref{2}et \ref{3}. Alors on a 

I) $\kq_{\zeta}$ est stable si et seulement si $\gamma -\zeta n > (2n+1)\lambda_{n}$.

II) Si $\kq_{\zeta}$ est stable, alors $\kn_{\zeta}$ est simple}.
 
 \begin{proof}
 Le th\'eor\`eme 2.7, \cite{bo-sp} donne (I).

Pour d\'emontrer II), supposons que $\kq_{\zeta}$ est stable. En tensorisant la suite \ref{3} par $\kn_{\zeta}$, on obtient la suite exacte suivante
$$0\longrightarrow \kn_{\zeta}\otimes \kn_{\zeta}^{*}\longrightarrow \kq_{\zeta}\otimes \kn_{\zeta}^{*} \stackrel{A}{\longrightarrow} \kn_{\zeta}^{*}(\gamma)\longrightarrow 0.$$

Alors on a
$$ h^{0}(\kn_{\zeta}\otimes \kn_{\zeta}^{*})\leq h^{0}(\kq_{\zeta}\otimes \kn_{\zeta}^{*}). $$

D'apr\`es le lemme \ref{0.2.3} et la proposition \ref{2.3}, on a que 
$$ 1\leq  h^{0}(\kn_{\zeta}\otimes \kn_{\zeta}^{*}).$$
 
On tensorise la suite duale de la suite \ref{3} par $\kq_{\zeta}$ ce qui permet d'obtenir la suite exacte suivante
$$ 0\longrightarrow \kq_{\zeta}(-\gamma)\stackrel{^{T}A}{\longrightarrow} \kq_{\zeta}\otimes \kq_{\zeta}^{*}\longrightarrow \kq_{\zeta}\otimes \kn_{\zeta}^{*}\longrightarrow 0.$$

Cependant, de la suite \ref{2} on a 
$$ h^{0}(\kq_{\zeta}(-\gamma))= h^{1}(\kq_{\zeta}(-\gamma))=0 .$$

Alors on a 
$$  h^{0}(\kn_{\zeta}^{*}\otimes \kq_{\zeta})=h^{0}(\kq_{\zeta}^{*}\otimes \kq_{\zeta}). $$

Mais, $ \kq_{\zeta}$ est stable. On en d\'eduit que $h^{0}(\kn_{\zeta}\otimes \kn_{\zeta}^{*})=1 $.

\end{proof}

\subsection{\bf Remarque}\label{2b}
Pour la suite exacte \ref{2}, on a la r\'esolution suivante
\begin{eqnarray} \label{4} 
0 \longrightarrow \ko_{\P^{2n+1}}(-(\gamma +\zeta) q)\longrightarrow 
\kh(-(\gamma +\zeta)(q-1))\longrightarrow \bigwedge ^{2}\kh(-(\gamma +\zeta)(q-2))\\
\nonumber \stackrel{a_{q-2}}{\longrightarrow}\bigwedge ^{3}\kh(-(\gamma +\zeta)(q-3)) 
\stackrel{a_{q-3}}{\longrightarrow} \ldots \\
\nonumber \ldots \stackrel{a_{3}}{\longrightarrow}\bigwedge 
^{q-2}\kh(-2(\gamma +\zeta)) \stackrel{a_{2}}{\longrightarrow}\bigwedge 
^{q-1}\kh(-(\gamma +\zeta)) \stackrel{a_{1}}{\longrightarrow}\bigwedge 
^{q}\kh\stackrel{a_{0}}{\longrightarrow}\bigwedge ^{q}\kq \longrightarrow  0
\end{eqnarray}

o\`u $1\leqslant q \leqslant 2n+1$. Pour la suite exacte \ref{2}, on a aussi la r\'esolution suivante 
\begin{eqnarray} \label{5}  
0 \longrightarrow \ko_{\P^{2n+1}}(-\gamma(2n+2-q)+\zeta (-(n+1)+q))\longrightarrow \bigwedge ^{2n+1}\kh_{\zeta}(-\gamma(2n+1-q)+\zeta q)\\
\nonumber \longrightarrow \bigwedge ^{2n}\kh_{\zeta}(-\gamma(2n-q)+\zeta q)\stackrel{a_{2n-1-q}}{\longrightarrow}\bigwedge ^{2n-1}\kh_{\zeta}(-\gamma(2n-1-q)+\zeta q) \stackrel{a_{2n-2-q}}{\longrightarrow} \ldots\\
\nonumber  \ldots \stackrel{a_{n+2-q}}{\longrightarrow}\bigwedge ^{n+2}\kh_{\zeta}(-\gamma(n+2-q)+\zeta q)\stackrel{a_{n+1-q}}{\longrightarrow}\bigwedge ^{n+1}\kh_{\zeta}(-\gamma(n+1-q) +\zeta q)\\
\nonumber \stackrel{a_{n-q}}{\longrightarrow} \bigwedge ^{n}\kh_{\zeta}(-\gamma(n-q)+\zeta q)\stackrel{a_{n-1-q}}{\longrightarrow} \ldots \\ 
\nonumber \ldots  \stackrel{a_{3}}{\longrightarrow}\bigwedge ^{q+3}\kh_{\zeta}(-3\gamma+\zeta q) \stackrel{a_{2}}{\longrightarrow}\bigwedge ^{q+2}\kh_{\zeta}(-2\gamma+\zeta q) \stackrel{a_{1}}{\longrightarrow}\bigwedge ^{q+1}\kh_{\zeta}(-\gamma+\zeta q)\stackrel{a_{0}}{\longrightarrow}\bigwedge ^{q}\kq_{\zeta}^{*} \longrightarrow  0
\end{eqnarray}

o\`u $1\leqslant q < 2n+1$. Pour la suite exacte duale de la suite \ref{3}
\begin{equation} \label{8}
0\longrightarrow \ko_{\P^{2n+1}}(-\gamma)\longrightarrow \kq_{\zeta}^{*}\longrightarrow  \kn_{\zeta}^{*}\longrightarrow 0,
\end{equation}

on a la r\'esolution suivante
\begin{eqnarray} \label{6} 
0 \longrightarrow \ko_{\P^{2n+1}}(- q\gamma)\longrightarrow 
\kq_{\zeta}^{*}(-\gamma(q-1))\longrightarrow \bigwedge ^{2}\kq_{\zeta}^{*}(-\gamma(q-2))\\
\nonumber  \stackrel{a_{q-2}}{\longrightarrow}\bigwedge ^{3}\kq_{\zeta}^{*}(-\gamma(q-3)) 
\stackrel{a_{q-3}}{\longrightarrow} \ldots\\
\nonumber  \ldots  \stackrel{a_{3}}{\longrightarrow}\bigwedge 
^{q-2}\kq_{\zeta}^{*}(-2\gamma) \stackrel{a_{2}}{\longrightarrow}\bigwedge 
^{q-1}\kq_{\zeta}^{*}(-\gamma) \stackrel{a_{1}}{\longrightarrow}\bigwedge 
^{q}\kq_{\zeta}^{*}\stackrel{a_{0}}{\longrightarrow}\bigwedge ^{q}\kn_{\zeta}^{*}
\longrightarrow  0
\end{eqnarray}

o\`u $1\leqslant q \leqslant 2n$.

 \newpage

\section{ \large\bf Stabilit\'e de fibr\'e de 0-corr\'elation pond\'er\'e}
\vspace{1cm}

Nous allons trouver les conditions n\'ecessaires et suffisantes pour que le fibr\'e de 0-corr\'elation pond\'er\'e $\kn_{\zeta}$ sur $\P^{2n+1}$ soit stable.

\subsection{\bf Lemme}\label{3.1} 
{\em Soient $p,\hspace{0.2 cm}q,\hspace{0.2 cm}k$ des entiers. Soit $\kh_{\zeta}$ un fibr\'e comme dans \ref{2s} pour $\zeta=0$.

I- Pour $k\geq \gamma$ et $1\leq q\leq 2n+1$, on a

$ \hspace{0.7 cm}$ $h^{0}(\bigwedge^{q}\kh_{\zeta}(-k))=0 $ si et seulement si  $ \gamma > \sum_{i=0}^{min(q-1,n)}\lambda_{n-i}  $.
En particulier, on a 
$$ h^{0}(\bigwedge^{q}\kh_{\zeta}(-k))=0 \hspace{0.2 cm} si \hspace{0.2 cm} \gamma > \sum_{i=0}^{n}\lambda_{i}  .$$

II- Pour $k\geq 2$ et $1\leq p,q\leq 2n+1$, on a
 
$ \hspace{0.7 cm}$ $h^{0}((\bigwedge^{p}\kh_{\zeta})\otimes \bigwedge^{q} \kh_{\zeta}(-k\gamma))=0 $ si  $ 
\gamma > \sum_{i=0}^{n}\lambda_{i}  $}.

\begin{proof}
 Pour d\'emontrer (I), on a que 
$$max\lbrace t \in \Z|\hspace{0.2 cm}  \ko_{\P^{2n+1}}(t)\subseteq \bigwedge^{q}\kh_{\zeta}\rbrace = \sum_{i=0}^{min(q-1,n)}\lambda_{n-i}.$$
 
Donc pour $k\geq \gamma$, on obtient 
$$h^{0}(\bigwedge^{q} \kh_{\zeta}(-k))=0 \hspace{0.2 cm} si \hspace{0.2 cm} et  \hspace{0.2 cm}seulement  \hspace{0.2 cm} si \hspace{0.2 cm} \gamma >\sum_{i=0}^{min(q-1,n)}\lambda_{n-i}. $$

En particulier, on a
$$h^{0}(\bigwedge^{q}\kh_{\zeta}(-k))=0 \hspace{0.2 cm} si \hspace{0.2 cm} \gamma > \sum_{i=0}^{n}\lambda_{i}  .$$

Pour d\'emontrer (II), on a que
$$max\lbrace t \in \Z|\hspace{0.2 cm}  \ko_{\P^{2n+1}}(t)\subseteq (\bigwedge^{p}\kh_{\zeta})\otimes (\bigwedge^{q} \kh_{\zeta})\rbrace = \sum_{i=0}^{min(q-1,n)}\lambda_{n-i}+ \sum_{i=0}^{min(p-1,n)}\lambda_{n-i}\leq  2\sum_{i=0}^{n}\lambda_{i}.$$
 
Donc pour $k\geq 2$, on obtient
$$h^{0}((\bigwedge^{p}\kh_{\zeta})\otimes \bigwedge^{q} \kh_{\zeta}(-k\gamma))=0 \hspace{0.2 cm} si \hspace{0.2 cm} \gamma >\sum_{i=0}^{n}\lambda_{i} . $$
 
\end{proof}

\subsection{\bf Lemme}\label{3.1.1} 
{\em Soient $p,\hspace{0.2 cm}q,\hspace{0.2 cm}k,\hspace{0.2 cm}b,\hspace{0.2 cm}a$ des entiers. Soit $\kh_{\zeta}$ un fibr\'e comme dans \ref{2s} pour $\zeta=1$.

I- Pour $k\geq \gamma$ et $1\leq q\leq 2n+1$ et $ 0\leq |a |  \leq n$, on a

$ \hspace{0.7 cm}$ $h^{0}(\bigwedge^{q}\kh_{\zeta}(-k+a))=0 $ si et seulement si  $ \gamma - n> \sum_{i=0}^{min(q-1,n)}\lambda_{n-i}  $.
En particulier, on a 
$$ h^{0}(\bigwedge^{q}\kh_{\zeta}(-k+a))=0 \hspace{0.2 cm} si \hspace{0.2 cm} \gamma -n > \sum_{i=0}^{n}\lambda_{i}  .$$

II- Pour $k\geq 2, 0\leq |b|\leq n$ et $1\leq p,q\leq 2n+1$, on a
 
$ \hspace{0.7 cm}$ $h^{0}((\bigwedge^{p}\kh_{\zeta})\otimes \bigwedge^{q} \kh_{\zeta}(-k\gamma-b))=0 $ si  $ 
\gamma - n> \sum_{i=0}^{n}\lambda_{i}  $}.

\begin{proof}
La d\'emonstration de ce lemme est tr\`es similaire \`a celle du lemme \ref{3.1}.

\end{proof}

\subsection{\bf Proposition}\label{3.2} 
{\em Soit $\kq_{\zeta}$ le fibr\'e de quotient pond\'er\'e sur $\P^{2n+1}$ qui est d\'efini par la suite exacte \ref{2} pour $\zeta=0$. Soient $1\leq q\leq n$, $0< i < 2n+1$ et $k\geq 0$ des entiers. Alors on a

I-  $h^{0}(\bigwedge^{q}\kq_{\zeta}^{*}(-m))=0$ si   
$\gamma > \sum_{v=0}^{n}\lambda_{v}$, $\forall m \in \N$.

II-  Pour $0\leq k\leq  q$, on a }
$$
h^{i}(\bigwedge^{q}\kq_{\zeta}^{*}(-k\gamma))=
\left\{
\begin{array}{ccc}
0 & : &  i\neq q\\
\epsilon_{k,q} & : & i=q\neq k\\
1  &  : & i=q= k
\end{array}
\right.
$$

{\em $ \hspace{0.7 cm}$ o\`u  $\epsilon_{k,q}=h^{1}\lbrace ker[\kh_{\zeta}(\gamma(q-1-k))\longrightarrow \ko_{\P^{2n+1}}(\gamma(q-k))]\rbrace$.

III-   Pour $ k> q$, on a  
$$h^{i}(\bigwedge^{q}\kq_{\zeta}^{*}(-k\gamma))=0.$$}

\begin{proof}
   
Pour d\'emontrer (I), on consid\`ere $1\leq q\leq n$, $0< i < 2n+1$ et $\forall m \in 
\N$ et $\zeta=0$. De la r\'esolution \ref{5}, on obtient la r\'esolution suivante
$$0 \longrightarrow \ko_{\P^{2n+1}}(-\gamma(2n+2-q)-m)\longrightarrow 
\kh_{\zeta}(-\gamma(2n+1-q)-m)\longrightarrow \bigwedge ^{2}\kh_{\zeta}(-\gamma(2n-q)-m)$$
$$\stackrel{a_{2n-1-q}}{\longrightarrow}\bigwedge ^{3}\kh_{\zeta}(-\gamma(2n-1-q)-m) 
\stackrel{a_{2n-2-q}}{\longrightarrow} \ldots$$
$$\ldots  \stackrel{a_{3}}{\longrightarrow}\bigwedge ^{q+3}\kh_{\zeta}(-3\gamma-m) 
\stackrel{a_{2}}{\longrightarrow}\bigwedge ^{q+2}\kh_{\zeta}(-2\gamma-m) 
\stackrel{a_{1}}{\longrightarrow}\bigwedge 
^{q+1}\kh_{\zeta}(-\gamma-m)\stackrel{a_{0}}{\longrightarrow}\bigwedge ^{q}\kq_{\zeta}^{*}(-m) 
\longrightarrow  0.$$

Soit $A_{r}=ker (a_{r})$ o\`u $r=0,1, \ldots ,2n-1-q$. Alors on obtient que 
$h^{j}(A_{0})=0$ pour $j\leq q$. Mais on a  
$$h^{0}(\bigwedge^{q+1} \kh_{\zeta}(-\gamma-m))=0 \hspace{0.2 cm} si \hspace{0.2 cm}  \gamma >\sum_{v=0}^{n}\lambda_{v} .$$

Donc on obtient, pour tout $m \in \N$, 
$$ h^{0}(\bigwedge^{q}\kq_{\zeta}^{*}(-m))=0\hspace{0.2 cm} si\hspace{0.2 cm}  \gamma > \sum_{v=0}^{n}\lambda_{v}.$$

De plus on a $h^{i}(\bigwedge^{q}\kq_{\zeta}^{*}(-m))=0 \hspace{0.2 cm} pour 
\hspace{0.2 cm} 0<i\leq q-1 $. De la r\'esolution duale de r\'esolution~\ref{4} on obtient la r\'esolution suivante, pour tout $q>0$,
$$0 \longrightarrow \bigwedge ^{q}\kq_{\zeta}^{*}(-k\gamma)\longrightarrow \bigwedge ^{q}\kh_{\zeta}(-k\gamma) 
\longrightarrow \bigwedge ^{q-1}\kh_{\zeta}(\gamma(1-k)) \stackrel{b_{q-2}}{\longrightarrow} 
\bigwedge ^{q-2}\kh_{\zeta}(\gamma(2-k))\stackrel{b_{q-3}}{\longrightarrow} \ldots$$
$$\ldots   \stackrel{b_{3}}{\longrightarrow}\bigwedge ^{3}\kh_{\zeta}(\gamma(q-3-k)) 
\stackrel{b_{2}}{\longrightarrow}\bigwedge ^{2}\kh_{\zeta}(\gamma(q-2-k))
\stackrel{b_{1}}{\longrightarrow}\kh_{\zeta}(\gamma(q-1-k))\stackrel{b_{0}}{\longrightarrow}\ko(\gamma(q-k))
 \longrightarrow 0.$$

On consid\`ere $B_{r}=ker (b_{r})$ pour $r=0,1, \ldots ,q-2$. Alors on obtient 
$$
h^{i}(B_{q-2})=
\left\{
\begin{array}{ccc}
0 & : & q\leq i\leq 2n\\
h^{1}(B_{0})=\epsilon_{k,q} & : & i=q-1\\
\end{array}
\right.
$$

Comme on a la suite exacte
$$0 \longrightarrow \bigwedge ^{q}\kq_{\zeta}^{*}(-k\gamma)\longrightarrow \bigwedge ^{q}\kh_{\zeta}(-k\gamma) 
\longrightarrow  B_{q-2} \stackrel{ }{\longrightarrow}0,$$

alors on a, pour $0\leq k\leq  q$ et $0< i < 2n+1$, 
$$
h^{i}(\bigwedge^{q}\kq_{\zeta}^{*}(-k\gamma))=
\left\{
\begin{array}{ccc}
0 & : & q\neq i\\
1 & : & i=q= k\\
\epsilon_{k,q}  &  : & sinon
\end{array}
\right.
$$

Pour $ k > q$ et $0< i < 2n+1$, on a  
$$h^{i}(\bigwedge^{q}\kq_{\zeta}^{*}(-k\gamma))=0.$$
  
\end{proof}

\subsection{\bf Proposition}\label{3.2.1} 
{\em Soit $\kq_{\zeta}$ le fibr\'e de quotient pond\'er\'e sur $\P^{2n+1}$ qui est d\'efini par la suite exacte \ref{2} pour $\zeta=1$. Soient $1\leq q\leq n$, $0< i < 2n+1$ et $k\geq 0$ et $0\leq  a \leq n$ des entiers. Alors on a

I- $h^{0}(\bigwedge^{q}\kq_{\zeta}^{*}(-k\gamma -a))=0$ si   
$\gamma-n > \sum_{v=0}^{n}\lambda_{v} $.

II- Pour $0\leq k,a\leq  q$, on a }
$$
h^{i}(\bigwedge^{q}\kq_{\zeta}^{*}(-k\gamma -a))=
\left\{
\begin{array}{ccc}
0 & : & q\neq i\\
1 & : & i=q= k=a\\
\epsilon_{k,q,a}  &  : & sinon
\end{array}
\right.
$$

{\em $ \hspace{0.7 cm}$ o\`u  $\epsilon_{k,q,a}=h^{1}\lbrace Ker[\kh_{\zeta}(\gamma(q-1-k)+q-1-a )\longrightarrow \ko_{\P^{2n+1}}(\gamma (q-k)+q-a)]\rbrace$.

III-  Pour $ k,a > q$, on a  
$$h^{i}(\bigwedge^{q}\kq_{\zeta}^{*}(-k\gamma -a))=0.$$}

\begin{proof}
La d\'emonstration de cette proposition est tr\`es similaire \`a celle de la proposition \ref{3.2}.

\end{proof}

\subsection{\bf Proposition}\label{3.3a} 
{\em Soit $\kn_{\zeta}$ le fibr\'e de \mbox{0-corr\'elation} pond\'er\'e sur $\P^{2n+1}$ qui est d\'efini par la suite exacte \ref{8} pour $\zeta=0$. Pour $1\leq 2j+1\leq n$ et $m\geq 0$ des entiers,
 
I- Si on a  $\hspace{0.2 cm}\gamma > \sum_{i=0}^{n}\lambda_{i}$, alors on obtient que $h^{0}(\bigwedge^{2j+1}\kn_{\zeta}(-m\gamma))=0$.

II- $h^{1}(\bigwedge^{2j+1}\kn_{\zeta}(-\gamma))=1$}.

\begin{proof}
    
Pour d\'emontrer (I), on fixe  $1\leq 2j+1\leq n$, $\forall  m \in \N$ et $\zeta=0$. On obtient, de la r\'esolution \ref{6}, la r\'esolution suivante
$$0 \longrightarrow \ko_{\P^{2n+1}}(- \gamma(2j+1+m))\longrightarrow 
\kq_{\zeta}^{*}(-\gamma(2j+m))\longrightarrow \bigwedge ^{2}\kq_{\zeta}^{*}(-\gamma(2j-1+m))$$
$$\stackrel{a_{2j-1}}{\longrightarrow}\bigwedge ^{3}\kq_{\zeta}^{*}(-\gamma(2j-2+m)) 
\stackrel{a_{2j-2}}{\longrightarrow} \ldots$$
$$\ldots   \stackrel{a_{2}}{\longrightarrow}\bigwedge 
^{2j}\kq_{\zeta}^{*}(-\gamma(1+m)) \stackrel{a_{1}}{\longrightarrow}\bigwedge 
^{2j+1}\kq_{\zeta}^{*}(-\gamma m)\stackrel{a_{0}}{\longrightarrow}\bigwedge 
^{2j+1}\kn_{\zeta}(-\gamma m) \longrightarrow    0.$$

On consid\`ere $A_{i}=ker (a_{i})$ pour $i=0,1, \ldots ,2j-1$. Les termes dans cette r\'esolution sont de la forme $\bigwedge ^{q}\kq_{\zeta}^{*}(-\gamma(k+m))$, o\`u $0\leq k\leq 2j$, $1\leq q \leq 2j+1$ tels que $k+q=2j+1$ et $k\neq q$. D'apr\`es la proposition \ref{3.2}, on obtient que $h^{k}(\bigwedge ^{2j+1-k}\kq_{\zeta}^{*}(-\gamma(k+m)))=0$. Comme 
$$h^{2j+1}(\ko_{\P^{2n+1}}(- \gamma(2j+1+m)))=0$$

alors on a $h^{i+1}(A_{i})=0 ,\hspace{0.2 cm} i=0,1, \ldots ,2j-1$. De la suite exacte 
$$0\longrightarrow A_{0}\longrightarrow\bigwedge ^{2j+1}\kq_{\zeta}^{*}(-\gamma m)\stackrel{a_{0}}{\longrightarrow}\bigwedge ^{2j+1}\kn_{\zeta}(-\gamma m) \longrightarrow    0$$ 

et comme $h^{0}(\bigwedge ^{2j+1}\kq_{\zeta}^{*}(-\gamma m))=0\hspace{0.2 cm}  si\hspace{0.2 cm} \gamma > \sum_{i=0}^{n}\lambda_{i}$, alors on obtient que
$$h^{0}(\bigwedge ^{2j+1}\kn_{\zeta}(-\gamma m))=0\hspace{0.2 cm} si\hspace{0.2 cm} \gamma > \sum_{i=0}^{n}\lambda_{i}.$$

Pour d\'emontrer (II), on fixe $  m=1$ dans la r\'esolution pr\'ec\'edente. On obtient la r\'esolution suivante
$$0 \longrightarrow \ko_{\P^{2n+1}}(- \gamma(2j+2))\longrightarrow 
\kq_{\zeta}^{*}(-\gamma(2j+1))\longrightarrow \bigwedge ^{2}\kq_{\zeta}^{*}(-\gamma(2j))$$
$$\stackrel{d_{2j-1}}{\longrightarrow}\bigwedge ^{3}\kq_{\zeta}^{*}(-\gamma(2j-1)) 
\stackrel{d_{2j-2}}{\longrightarrow} \ldots $$
$$ \ldots \stackrel{d_{j+2}}{\longrightarrow}\bigwedge^{j} \kq_{\zeta}^{*}(- 
\gamma(j+2))\stackrel{d_{j+1}}{\longrightarrow} 
\bigwedge^{j+1}\kq_{\zeta}^{*}(-\gamma(j+1))\stackrel{d_{j}}{\longrightarrow} \bigwedge 
^{j+2}\kq_{\zeta}^{*}(-\gamma(j))\stackrel{d_{j-1}}{\longrightarrow} \ldots$$
$$ \ldots  \stackrel{d_{3}}{\longrightarrow}\bigwedge 
^{2j-1}\kq_{\zeta}^{*}(-3\gamma) \stackrel{d_{2}}{\longrightarrow}\bigwedge 
^{2j}\kq_{\zeta}^{*}(-2\gamma) \stackrel{d_{1}}{\longrightarrow}\bigwedge 
^{2j+1}\kq_{\zeta}^{*}(-\gamma)\stackrel{d_{0}}{\longrightarrow}\bigwedge 
^{2j+1}\kn_{\zeta}(-\gamma) \longrightarrow    0.$$

On consid\`ere $D_{\alpha}=ker (d_{\alpha})$ pour $\alpha=0,1, \ldots ,2j-1$. Les termes dans cette r\'esolution sont de la forme $\bigwedge ^{q}\kq_{\zeta}^{*}(-k\gamma)$ o\`u $1\leq k, q \leq 2j+1$ tels que  $k+q=2j+2$.

On a deux cas: 

Premier cas, si $k\geq q$. Dans ce cas on obtient, d'apr\`es la proposition \ref{3.2}, que
$$h^{i}(\bigwedge ^{q}\kq_{\zeta}^{*}(-k\gamma))=0 \hspace{0.2 cm}  pour \hspace{0.2 cm} i\neq 0,\hspace{0.2 cm} 2n+1.$$

Donc \c{c}a nous donne que $h^{i}(D_{\alpha})=0  \hspace{0.2 cm} pour  \hspace{0.2 cm} \alpha \geq j, \hspace{0.2 cm} i\leq 2n-j $, sauf pour $ i=\alpha $. On a   
$$h^{j+1}(\bigwedge ^{j+1}\kq_{\zeta}^{*}(-\gamma(j+1)))=1 \hspace{0.2 cm} , \hspace{0.2 cm} h^{i}(D_{j})=0   \hspace{0.2 cm} , \hspace{0.2 cm} i=j+1,j+2.$$

De la suite exacte suivante
$$0\longrightarrow D_{j}\longrightarrow\bigwedge ^{j+1}\kq_{\zeta}^{*}(-\gamma (j+1))\longrightarrow  D_{j-1} \longrightarrow    0,$$

on obtient que $h^{j+1}(D_{j-1})=1  $. Donc on obtient $h^{\alpha +2}(D_{\alpha})=1\hspace{0.2 cm}  pour \hspace{0.2 cm}\alpha =0,1, \ldots ,j-1$.

Deuxi\`eme cas, si $k<q$. D'apr\`es la proposition \ref{3.2}, on obtient $h^{i}(\bigwedge ^{q}\kq_{\zeta}^{*}(-k\gamma))=0 $  pour $ i=k,\hspace{0.2 cm}k-1\hspace{0.2 cm},1\leq k\leq j+1 $. Alors on obtient que $h^{\alpha +2}(D_{\alpha})=1\hspace{0.2 cm}  pour \hspace{0.2 cm}\alpha =0,1, \ldots ,j-1$. 

Dans les deux cas: d'apr\`es la proposition \ref{3.2}, on obtient que 
$$h^{i}(\bigwedge ^{2j+1}\kq_{\zeta}^{*}(-\gamma))=0 \hspace{0.2 cm}  pour \hspace{0.2 cm}i=1,2.$$

Donc d'apr\`es la suite exacte suivante
$$0\longrightarrow D_{0}\longrightarrow\bigwedge ^{2j+1}\kq_{\zeta}^{*}(-\gamma )\stackrel{a_{0}}{\longrightarrow}\bigwedge ^{2j+1}\kn_{\zeta}(-\gamma ) \longrightarrow 0,$$

on obtient que $h^{1}(\bigwedge ^{2j+1} \kn_{\zeta}(- \gamma))=1$.

\end{proof}

\subsection{\bf Proposition}\label{3.3a.1}
{\em Soit $\kn_{\zeta}$ le fibr\'e de \mbox{0-corr\'elation} pond\'er\'e sur $\P^{2n+1}$ qui est d\'efini par la suite exacte \ref{8} pour $\zeta=1$. Pour $1\leq 2j+1\leq n$, $m\geq 0$ et $0\leq  a \leq n$ des entiers, 
 
I- Si on a $\hspace{0.2 cm}\gamma-n > \sum_{i=0}^{n}\lambda_{i}$, alors on obtient que $h^{0}(\bigwedge^{2j+1}\kn_{\zeta}^{*}(-m\gamma-a ))=0$.

II- $h^{1}(\bigwedge^{2j+1}\kn_{\zeta}^{*}(-\gamma-(j+1)))=1$.}

\begin{proof}
    
Pour d\'emontrer (I), on fixe  $1\leq 2j+1\leq n$ et $\zeta=1$. On obtient, de la r\'esolution \ref{6}, la r\'esolution suivante
$$0 \longrightarrow \ko_{\P^{2n+1}}(- \gamma(2j+1+m)-\alpha)\longrightarrow \kq_{\zeta}^{*}(- \gamma(2j+m)-\alpha)\longrightarrow \bigwedge ^{2}\kq_{\zeta}^{*}(- \gamma(2j-1+m)-\alpha)$$
$$\stackrel{a_{2j-1}}{\longrightarrow}\bigwedge ^{3}\kq_{\zeta}^{*}(- \gamma(2j-2+m)-\alpha) \stackrel{a_{2j-2}}{\longrightarrow}\ldots $$
$$\stackrel{a_{j+1}}{\longrightarrow}\bigwedge ^{j+1}\kq_{\zeta}^{*}(- \gamma(j+m)-\alpha) \stackrel{a_{j}}{\longrightarrow} \bigwedge ^{j+2}\kq_{\zeta}^{*}(- \gamma(j-1+m)-\alpha) \stackrel{a_{j-1}}{\longrightarrow}\ldots$$
$$\ldots  \stackrel{a_{2}}{\longrightarrow}\bigwedge ^{2j}\kq_{\zeta}^{*}(- \gamma(1+m)-\alpha) \stackrel{a_{1}}{\longrightarrow}\bigwedge ^{2j+1}\kq_{\zeta}^{*}(- \gamma m-\alpha)\stackrel{a_{0}}{\longrightarrow}\bigwedge ^{2j+1}\kn_{\zeta}^{*}(- \gamma m-\alpha) \longrightarrow    0$$

On consid\`ere  $A_{i}=Ker (a_{i})$ pour $i=0,1,\ldots ,2j-1$. Les termes dans cette r\'esolution sont de la forme $\bigwedge ^{q}\kq_{\zeta}^{*}(-\gamma(k+m)-\alpha)$, o\`u $0\leq k\leq 2j$, $1\leq q \leq 2j+1$ tels que  $k+q=2j+1$ et $k\neq q$. D'apr\`es la proposition \ref{3.2.1}, on obtient que $h^{k}(\bigwedge ^{2j+1-k}\kq_{\zeta}^{*}(-\gamma(k+m)-\alpha))=0$. Comme 
$$h^{2j+1}(\ko_{\P^{2n+1}}(- \gamma(2j+1+m)-\alpha))=0$$

alors on a $h^{i+1}(A_{i})=0 \hspace{0.2 cm},\hspace{0.2 cm} i=0,1,..,2j-1$. De la suite exacte 
$$0\longrightarrow A_{0}\longrightarrow\bigwedge ^{2j+1}\kq_{\zeta}^{*}(-\gamma m-\alpha)\stackrel{a_{0}}{\longrightarrow}\bigwedge ^{2j+1}\kn_{\zeta}^{*}(-\gamma m-\alpha) \longrightarrow    0$$ 

et comme $h^{0}(\bigwedge ^{2j+1}\kq_{\zeta}^{*}(-\gamma m))=0\hspace{0.2 cm} ,\hspace{0.2 cm} si\hspace{0.2 cm} \gamma -n> \sum_{i=0}^{n}\lambda_{i}$, alors on obtient que

$$h^{0}(\bigwedge ^{2j+1}\kn_{\zeta}^{*}(-\gamma m-\alpha))=0\hspace{0.2 cm} ,\hspace{0.2 cm} si\hspace{0.2 cm} \gamma-n > \sum_{i=0}^{n}\lambda_{i}.$$

Pour d\'emontrer (II), on fixe $  m=1, \alpha =j+1$ dans la r\'esolution pr\'ec\'edente; on obtient la r\'esolution suivante
$$0 \longrightarrow \ko_{\P^{2n+1}}(- \gamma(2j+2)-(j+1))\longrightarrow \kq_{\zeta}^{*}(-\gamma(2j+1)-(j+1))\longrightarrow \bigwedge ^{2}\kq_{\zeta}^{*}(-\gamma(2j)-(j+1))$$
$$\stackrel{d_{2j-1}}{\longrightarrow}\bigwedge ^{3}\kq_{\zeta}^{*}(-\gamma(2j-1)-(j+1)) \stackrel{d_{2j-2}}{\longrightarrow}\ldots$$
$$ \ldots \stackrel{d_{j+2}}{\longrightarrow}\bigwedge^{j} \kq_{\zeta}^{*}(- \gamma(j+2)-(j+1))\stackrel{d_{j+1}}{\longrightarrow} \bigwedge^{j+1}\kq_{\zeta}^{*}(-\gamma(j+1)-(j+1))\stackrel{d_{j}}{\longrightarrow} \bigwedge ^{j+2}\kq_{\zeta}^{*}(-\gamma(j)-(j+1))\stackrel{d_{j-1}}{\longrightarrow}\ldots$$
$$\ldots \stackrel{d_{2}}{\longrightarrow}\bigwedge ^{2j}\kq_{\zeta}^{*}(-2\gamma-(j+1)) \stackrel{d_{1}}{\longrightarrow}\bigwedge ^{2j+1}\kq_{\zeta}^{*}(-\gamma-(j+1))\stackrel{d_{0}}{\longrightarrow}\bigwedge ^{2j+1}\kn_{\zeta}^{*}(-\gamma-(j+1)) \longrightarrow    0.$$

On consid\`ere  $D_{i}=Ker (d_{i})$ pour $i=0,1,\ldots ,2j-1$. Les termes dans cette r\'esolution sont de

la forme $\bigwedge ^{q}\kq_{\zeta}^{*}(-k\gamma-(j+1))$, o\`u $1\leq k, q \leq 2j+1$ tels que  $k+q=2j+2$. 

On a deux cas:

Premier cas, si $k\geq q$. Dans ce cas, on obtient, d'apr\`es la proposition \ref{3.2.1}, que
$$h^{i}(\bigwedge ^{q}\kq_{\zeta}^{*}(-k\gamma)-(j+1))=0 \hspace{0.2 cm}  pour \hspace{0.2 cm} i\neq 0,\hspace{0.2 cm} 2n+1.$$

Donc \c{c}a nous donne que $h^{i}(D_{\beta})=0  \hspace{0.2 cm} pour  \hspace{0.2 cm} \beta \geq j, \hspace{0.2 cm} i\leq 2n-j, \hspace{0.2 cm} sauf \hspace{0.2 cm} i=\beta $. On a 
$$h^{j+1}(\bigwedge ^{j+1}\kq_{\zeta}^{*}(-\gamma(j+1)-(j+1)))=1 \hspace{0.2 cm} , \hspace{0.2 cm} h^{i}(D_{j})=0   \hspace{0.2 cm} , \hspace{0.2 cm} i=j+1,j+2$$

et on a la suite exacte suivante
$$0\longrightarrow D_{j}\longrightarrow\bigwedge ^{j+1}\kq_{\zeta}^{*}(-\gamma (j+1)-(j+1))\longrightarrow  D_{j-1} \longrightarrow    0$$

On obtient que $h^{j+1}(D_{j-1})=1  $. Donc on obtient $h^{\beta +2}(D_{\beta})=1\hspace{0.2 cm}  pour \hspace{0.2 cm}\beta =0,1,\ldots ,j-1$. 

Deuxi\`eme cas, $k<q$. D'apr\`es la proposition \ref{3.2.1}, on obtient 
$$h^{i}(\bigwedge ^{q}\kq_{\zeta}^{*}(-k\gamma)-(j+1))=0 \hspace{0.2 cm}  pour \hspace{0.2 cm} i=k,\hspace{0.2 cm}k-1\hspace{0.2 cm},1\leq k\leq j+1 .$$

Alors on obtient que $h^{\beta +2}(D_{\beta})=1\hspace{0.2 cm}  pour \hspace{0.2 cm}\beta =0,1,\ldots ,j-1$. Dans les deux cas: d'apr\`es la proposition \ref{3.2.1}, on obtient que 
$$h^{i}(\bigwedge ^{2j+1}\kq_{\zeta}^{*}(-\gamma-(j+1)))=0 \hspace{0.2 cm}  pour \hspace{0.2 cm}i=1,2.$$

Donc d'apr\`es la suite exacte suivante
$$0\longrightarrow D_{0}\longrightarrow\bigwedge ^{2j+1}\kq_{\zeta}^{*}(-\gamma-(j+1) )\stackrel{a_{0}}{\longrightarrow}\bigwedge ^{2j+1}\kn_{\zeta}^{*}(-\gamma-(j+1) ) \longrightarrow 0,$$

on obtient que $h^{1}(\bigwedge ^{2j+1} \kn_{\zeta}^{*}(- \gamma -(j+1)))=1$.
 
\end{proof}

\subsection{\bf Proposition}\label{3.4} 
{\em Soit $\kq_{\zeta}$ le fibr\'e de quotient pond\'er\'e sur $\P^{2n+1}$ qui est d\'efini par la suite exacte \ref{2} pour $\zeta=0$. Soient $ \alpha\in \Z $, $ 1\leq k\leq n+1$ un entier. Si on a  $\gamma > \sum_{i=0}^{n}\lambda_{i}$, alors on obtient  

- $Pour\hspace{0.2 cm} \alpha =1$,
$$ h^{i}(\kq_{\zeta}^{*}\otimes \bigwedge^{k}\kh_{\zeta}(-\alpha \gamma))=0, \text{ o\`u } \hspace{0.1 cm} 2\leq \hspace{0.2 cm}i\leq   \hspace{0.2 cm} 2n, \hspace{0.2 cm} ou \hspace{0.2 cm} bien\hspace{0.2 cm}  i=0 .$$

- $Pour\hspace{0.2 cm} \alpha < 1$, 
$$  h^{i}(\kq_{\zeta}^{*}\otimes \bigwedge^{k}\kh_{\zeta}(-\alpha \gamma))=0,\text{ o\`u } \hspace{0.1 cm} 2\leq \hspace{0.2 cm}i\leq   \hspace{0.2 cm} 2n. $$

- $Pour\hspace{0.2 cm} \alpha >1$, 
$$  h^{i}(\kq_{\zeta}^{*}\otimes \bigwedge^{k}\kh_{\zeta}(-\alpha \gamma))=0,\text{ o\`u } \hspace{0.1 cm} 0\leq \hspace{0.2 cm}i\leq   \hspace{0.2 cm} 2n. $$}
 
\begin{proof}
    
On fixe $\gamma > \sum_{i=0}^{n}\lambda_{i}$. Soient $ 1\leq k\leq n+1 \hspace{0.2 cm}un\hspace{0.2 cm} entier\hspace{0.2 cm} et\hspace{0.2 cm}\alpha\in \Z $ et $\zeta=0$. On~veut calculer les groupes cohomologiques du fibr\'e 
$$\kq_{\zeta}^{*}\otimes \bigwedge^{k}\kh_{\zeta}(-\alpha \gamma) .$$

De la suite exacte suivante 
$$0\longrightarrow \kq_{\zeta}^{*} \otimes \bigwedge^{k}\kh_{\zeta}(-\alpha \gamma) \longrightarrow \kh_{\zeta} \otimes \bigwedge^{k}\kh_{\zeta}(-\alpha \gamma) \longrightarrow  \bigwedge^{k}\kh_{\zeta}(\gamma (-\alpha +1) ) \longrightarrow 0$$

et d'apr\`es le lemme \ref{3.1}: 

$si\hspace{0.2 cm} \alpha \leq 1$, on a
$$ h^{i}(\kq_{\zeta}^{*}\otimes \bigwedge^{k}\kh_{\zeta}(-\alpha \gamma))=0,\hspace{0.2 cm} 2\leq \hspace{0.2 cm}i\leq  \hspace{0.2 cm} 2n, $$

$si\hspace{0.2 cm} \alpha > 1$, on a
$$  h^{i}(\kq_{\zeta}^{*}\otimes \bigwedge^{k}\kh_{\zeta}(-\alpha \gamma))=0,\hspace{0.2 cm} 1\leq \hspace{0.2 cm}i\leq   \hspace{0.2 cm} 2n. $$

Pour calculer $h^{0}(\kq_{\zeta}^{*}\otimes \bigwedge^{k}\kh_{\zeta}(-\alpha \gamma))$, on va 
utiliser la r\'esolution \ref{5} avec $q=1$. Donc on obtient la r\'esolution 
suivante
$$0 \longrightarrow \bigwedge^{k}\kh_{\zeta}(-\gamma(\alpha +2n+1) )\longrightarrow 
\kh_{\zeta}\otimes \bigwedge^{k}\kh_{\zeta}(-\gamma(\alpha +2n) )\longrightarrow$$
$$(\bigwedge ^{2}\kh_{\zeta})\otimes \bigwedge^{k}\kh_{\zeta}(-\gamma(\alpha +2n-1) ) 
 \stackrel{d_{2n-2}}{\longrightarrow}(\bigwedge ^{3}\kh_{\zeta})\otimes 
 \bigwedge^{k}\kh_{\zeta}(-\gamma(\alpha +2n-2) ) 
 \stackrel{d_{2n-3}}{\longrightarrow} \ldots $$
$$ \ldots  \stackrel{d_{2}}{\longrightarrow}(\bigwedge ^{3}\kh_{\zeta})\otimes 
\bigwedge^{k}\kh_{\zeta}(-\gamma(\alpha +2) ) 
\stackrel{d_{1}}{\longrightarrow}(\bigwedge ^{2}\kh_{\zeta})\otimes 
\bigwedge^{k}\kh_{\zeta}(-\gamma(\alpha +1) 
)\stackrel{d_{0}}{\longrightarrow}\kq_{\zeta}^{*}\otimes \bigwedge^{k}\kh_{\zeta}(-\alpha 
\gamma) \longrightarrow  0.$$

On consid\`ere  $D_{j}=ker (d_{j})$ pour $j=0,1, \ldots ,2n-2$. D'apr\`es le lemme \ref{3.1}, on obtient que $h^{1}(D_{0})=0 $ pour $\alpha \in \Z$. D'apr\`es le lemme \ref{3.1} et de la suite exacte suivante
$$0\longrightarrow D_{0} \longrightarrow (\bigwedge ^{2}\kh_{\zeta})\otimes \bigwedge^{k}\kh_{\zeta}(-\gamma(\alpha +1) )\stackrel{d_{0}}{\longrightarrow}\kq_{\zeta}^{*}\otimes \bigwedge^{k}\kh_{\zeta}(-\alpha \gamma) \longrightarrow  0$$

on obtient que, $si\hspace{0.2 cm} \alpha \geq1$,  
$$ h^{0}(\kq_{\zeta}^{*}\otimes \bigwedge^{k}\kh_{\zeta}(-\alpha \gamma))=0 .$$

\end{proof}

\subsection{\bf Proposition}\label{3.4.1}
{\em Soit $\kq_{\zeta}$ le fibr\'e de quotient pond\'er\'e sur $\P^{2n+1}$ qui est d\'efini par la suite exacte \ref{2} pour $\zeta=1$. Soient $ \alpha\in \Z $, $ 2\leq k\leq 2n+1$ et $|b|<n$ des entiers. Si on a  $\gamma -n> \sum_{i=0}^{n}\lambda_{i}$, alors on obtient  

- $Pour\hspace{0.2 cm} \alpha =1$, 
$$ h^{i}(\kq_{\zeta}^{*}\otimes \bigwedge^{k}\kh_{\zeta}(-\alpha \gamma-b))=0, \text{ o\`u }\hspace{0.1 cm} 2\leq \hspace{0.2 cm}i\leq   \hspace{0.2 cm} 2n,\hspace{0.2 cm} ou \hspace{0.2 cm} bien \hspace{0.2 cm} i=0 $$

- $Pour\hspace{0.2 cm} \alpha < 1$, 
$$  h^{i}(\kq_{\zeta}^{*}\otimes \bigwedge^{k}\kh_{\zeta}(-\alpha \gamma-b))=0,\text{ o\`u }\hspace{0.1 cm} 2\leq \hspace{0.2 cm}i\leq   \hspace{0.2 cm} 2n. $$

- $Pour\hspace{0.2 cm} \alpha >1$, 
$$  h^{i}(\kq_{\zeta}^{*}\otimes \bigwedge^{k}\kh_{\zeta}(-\alpha \gamma-b))=0,\text{ o\`u }\hspace{0.1 cm} 0\leq \hspace{0.2 cm}i\leq   \hspace{0.2 cm} 2n. $$}
 
\begin{proof}
La d\'emonstration de cette proposition est tr\`es similaire \`a celle de la proposition \ref{3.4}.
\end{proof}

\subsection{\bf Proposition}\label{3.5a} 
{\em Soit $\kq_{\zeta}$ le fibr\'e de quotient pond\'er\'e sur $\P^{2n+1}$ qui est d\'efini par la suite exacte \ref{2} pour $\zeta=0$. Soient $ \alpha\geq 0 \hspace{0.2 cm} et \hspace{0.2 cm} 1 \leq q \leq n $ des entiers. Si on a  $\gamma > \sum_{i=0}^{n}\lambda_{i}$, alors on obtient

- Pour $\alpha =0$,
$$ h^{i}(\kq_{\zeta}^{*}\otimes\bigwedge ^{q}\kq_{\zeta}^{*}(-\alpha \gamma ))=0, \text{ o\`u } \hspace{0.1 cm} 2\leq \hspace{0.2 cm}i\leq   \hspace{0.2 cm} q-1, \text{ ou bien }  i=0 \hspace{0.1 cm} \text{ ou bien } q+1\leq \hspace{0.2 cm}i\leq   \hspace{0.2 cm} 2n. $$

- Pour $\alpha >0$,
$$  h^{i}(\kq_{\zeta}^{*}\otimes\bigwedge ^{q}\kq_{\zeta}^{*}(-\alpha \gamma ))=0, \text{ o\`u } \hspace{0.1 cm} 0\leq \hspace{0.2 cm}i\leq   \hspace{0.2 cm} q-1 \text{ ou bien } q+1\leq \hspace{0.2 cm}i\leq   \hspace{0.2 cm} 2n. $$}
 
\begin{proof}

On fixe $\gamma > \sum_{i=0}^{n}\lambda_{i}$. Soient $ \alpha\geq 0 \hspace{0.2 cm} et \hspace{0.2 cm} 1 \leq q \leq n $ des entiers et $\zeta=0$. On veut calculer les groupes cohomologiques du fibr\'e  
$$\kq_{\zeta}^{*}\otimes\bigwedge ^{q}\kq_{\zeta}^{*}(-\alpha \gamma ).$$

De la r\'esolution \ref{5}, on obtient la r\'esolution suivante
$$0 \longrightarrow \kq_{\zeta}^{*}(-\gamma(\alpha +2n+2-q))\longrightarrow 
\kq_{\zeta}^{*}\otimes\kh_{\zeta}(-\gamma(\alpha +2n+1-q))\longrightarrow 
\kq_{\zeta}^{*}\otimes\bigwedge ^{2}\kh_{\zeta}(-\gamma(\alpha +2n-q))$$
$$\stackrel{b_{2n-1-q}}{\longrightarrow}  \kq_{\zeta}^{*}\otimes\bigwedge 
^{3}\kh_{\zeta}(-\gamma(\alpha +2n-1-q)) \stackrel{b_{2n-2-q}}{\longrightarrow} \ldots $$
$$ \ldots  \stackrel{b_{2}}{\longrightarrow}\kq_{\zeta}^{*}\otimes\bigwedge 
^{q+2}\kh_{\zeta}(-\gamma(\alpha +2)) 
\stackrel{b_{1}}{\longrightarrow}\kq_{\zeta}^{*}\otimes\bigwedge 
^{q+1}\kh_{\zeta}(-\gamma(\alpha 
+1))\stackrel{b_{0}}{\longrightarrow}\kq_{\zeta}^{*}\otimes\bigwedge 
^{q}\kq_{\zeta}^{*}(-\alpha \gamma ) \longrightarrow  0.$$

On consid\`ere  $B_{m}=ker (b_{m})$ pour $m=0,1, \ldots ,2n-1-q$. D'apr\`es la proposition \ref{3.4}, on a alors 
$$h^{i}(B_{0})=0,\text{ o\`u } \hspace{0.1 cm} 0\leq i\leq q.$$

De la suite exacte suivante 
$$0\longrightarrow B_{0}\longrightarrow \kq_{\zeta}^{*}\otimes\bigwedge ^{q+1}\kh_{\zeta}(-\gamma(\alpha +1))\stackrel{b_{0}}{\longrightarrow}\kq_{\zeta}^{*}\otimes\bigwedge ^{q}\kq_{\zeta}^{*}(-\alpha \gamma ) \longrightarrow  0$$

et d'apr\`es la proposition \ref{3.4} on obtient, si $  \alpha =0$,
$$ h^{i}(\kq_{\zeta}^{*}\otimes\bigwedge ^{q}\kq_{\zeta}^{*}(-\alpha \gamma ))=0,\text{ o\`u } \hspace{0.1 cm} 2\leq \hspace{0.2 cm}i\leq   \hspace{0.2 cm} q-1,\text{ ou bien } i=0.$$

On obtient aussi, si $\alpha >0$,
$$  h^{i}(\kq_{\zeta}^{*}\otimes\bigwedge ^{q}\kq_{\zeta}^{*}(-\alpha \gamma ))=0,\text{ o\`u } \hspace{0.1 cm} 0\leq \hspace{0.2 cm}i\leq   \hspace{0.2 cm} q-1. $$

On a la r\'esolution suivante
$$0 \longrightarrow\kq_{\zeta}^{*}\otimes \bigwedge ^{q}\kq_{\zeta}^{*}(-\alpha 
\gamma)\longrightarrow \kq_{\zeta}^{*}\otimes \bigwedge ^{q}\kh_{\zeta}(-\alpha \gamma) 
\longrightarrow \kq_{\zeta}^{*}\otimes \bigwedge ^{q-1}\kh_{\zeta}(\gamma(-\alpha + 1) 
)\stackrel{a_{q-2}}{\longrightarrow}$$
$$\kq_{\zeta}^{*}\otimes \bigwedge ^{q-2}\kh_{\zeta}(\gamma(-\alpha + 2) 
)\stackrel{a_{q-3}}{\longrightarrow} \ldots$$
$$ \ldots \stackrel{a_{2}}{\longrightarrow}\kq_{\zeta}^{*}\otimes \bigwedge 
^{2}\kh_{\zeta}(\gamma(-\alpha + q-2) ) \stackrel{a_{1}}{\longrightarrow}\kq_{\zeta}^{*}\otimes 
\kh_{\zeta}(\gamma(-\alpha + q-1) 
)\stackrel{a_{0}}{\longrightarrow}\kq_{\zeta}^{*}(\gamma(-\alpha + q) ) \longrightarrow 
0.$$

On consid\`ere  $A_{r}=ker (a_{r})$ pour $r=0,1, \ldots ,q-2$. D'apr\`es la proposition \ref{3.4}, on a alors 
$$h^{i}(A_{q-2})=0,\text{ o\`u } \hspace{0.1 cm} q\leq i\leq 2n.$$

De la suite exacte suivante 
$$0 \longrightarrow\kq_{\zeta}^{*}\otimes \bigwedge ^{q}\kq_{\zeta}^{*}(-\alpha \gamma)\longrightarrow \kq_{\zeta}^{*}\otimes \bigwedge ^{q}\kh_{\zeta}(-\alpha \gamma)\longrightarrow A_{q-2}\longrightarrow  0$$
  
et d'apr\`es la proposition \ref{3.4} on obtient, si $ \alpha\geq 0 $,
$$ h^{i}(\kq_{\zeta}^{*}\otimes\bigwedge ^{q}\kq_{\zeta}^{*}(-\alpha \gamma ))=0, \text{ o\`u } \hspace{0.1 cm} q+1\leq \hspace{0.2 cm}i\leq   \hspace{0.2 cm} 2n. $$

\end{proof}

\subsection{\bf Proposition}\label{3.5a.1} 
{\em Soit $\kq_{\zeta}$ le fibr\'e de quotient pond\'er\'e sur $\P^{2n+1}$ qui est d\'efini par la suite exacte \ref{2} pour $\zeta=1$. Soient $ \alpha\geq 0 \hspace{0.2 cm} et \hspace{0.2 cm} 1 \leq q \leq n $ et $2\leq d\leq n $ des entiers. Si on a  $\gamma-n > \sum_{i=0}^{n}\lambda_{i}$, alors on obtient

- Pour $\alpha =0$,
$$ h^{i}(\kq_{\zeta}^{*}\otimes\bigwedge ^{q}\kq_{\zeta}^{*}(-\alpha \gamma-d ))=0,\text{ o\`u } \hspace{0.1 cm} 2\leq \hspace{0.2 cm}i\leq   \hspace{0.2 cm} q-1,\hspace{0.2 cm} ou \hspace{0.2 cm} i=0 \hspace{0.2 cm} ou \hspace{0.2 cm} q+1\leq \hspace{0.2 cm}i\leq   \hspace{0.2 cm} 2n $$

- Pour $\alpha >0$,
$$  h^{i}(\kq_{\zeta}^{*}\otimes\bigwedge ^{q}\kq_{\zeta}^{*}(-\alpha \gamma -d))=0,\text{ o\`u } \hspace{0.1 cm} 0\leq \hspace{0.2 cm}i\leq   \hspace{0.2 cm} q-1 \hspace{0.2 cm} ou \hspace{0.2 cm} q+1\leq \hspace{0.2 cm}i\leq   \hspace{0.2 cm} 2n. $$}
 
\begin{proof}
La d\'emonstration de cette proposition est tr\`es similaire \`a celle de la proposition \ref{3.5a}.
\end{proof}

\subsection{\bf Proposition}\label{3.6} 
{\em Soit $\kn_{\zeta}$ le fibr\'e de \mbox{0-corr\'elation} pond\'er\'e sur $\P^{2n+1}$ qui est d\'efini par la suite exacte \ref{8} pour $\zeta=0$. Soit $1\leq 2j+1\leq n$. Si on a $\gamma > \sum_{i=0}^{n}\lambda_{i}$, alors on obtient que 
$$h^{0}(\kn_{\zeta}\otimes \bigwedge^{2j+1}\kn_{\zeta})=1.$$}
 
\begin{proof}
On fixe $1\leq 2j+1\leq n$ et $\gamma > \sum_{i=0}^{n}\lambda_{i}$ et $\zeta=0$. D'apr\`es le lemme \ref{0.2.3} et la proposition \ref{2.3}, on obtient que 
$$h^{0}(\kn_{\zeta}\otimes \bigwedge^{2j+1}\kn_{\zeta})\geq 1.$$ 

De la suite exacte suivante
$$0\longrightarrow \bigwedge^{2j+1}\kn_{\zeta}(-\gamma)\longrightarrow \kq_{\zeta}^{*}\otimes \bigwedge^{2j+1}\kn_{\zeta}\longrightarrow  \kn_{\zeta}\otimes \bigwedge^{2j+1}\kn_{\zeta}\longrightarrow 0$$

et suivant la proposition \ref{3.3a}, il suffit de d\'emontrer que 
$$h^{0}(\kq_{\zeta}^{*}\otimes \bigwedge^{2j+1}\kn_{\zeta})=0$$

pour que l'on ait
$$h^{0}(\kn_{\zeta}\otimes \bigwedge^{2j+1}\kn_{\zeta})=1.$$

De la r\'esolution \ref{6}, on obtient la r\'esolution suivante
$$0 \longrightarrow \kq_{\zeta}^{*}(- \gamma(2j+1))\longrightarrow \kq_{\zeta}^{*}\otimes 
\kq_{\zeta}^{*}(-\gamma(2j))\longrightarrow \kq_{\zeta}^{*}\otimes \bigwedge 
^{2}\kq_{\zeta}^{*}(-\gamma(2j-1)) \stackrel{\mu_{2j-1}}{\longrightarrow}$$
$$\kq_{\zeta}^{*}\otimes \bigwedge ^{3}\kq_{\zeta}^{*}(-\gamma(2j-2)) 
\stackrel{\mu_{2j-2}}{\longrightarrow} \ldots $$
$$ \ldots \stackrel{\mu_{j+2}}{\longrightarrow} \kq_{\zeta}^{*}\otimes \bigwedge^{j} 
\kq_{\zeta}^{*}(- \gamma(j+1))\stackrel{\mu_{j+1}}{\longrightarrow} \kq_{\zeta}^{*}\otimes 
\bigwedge^{j+1}\kq_{\zeta}^{*}(-\gamma(j))\stackrel{\mu_{j}}{\longrightarrow} 
\kq_{\zeta}^{*}\otimes \bigwedge 
^{j+2}\kq_{\zeta}^{*}(-\gamma(j-1))\stackrel{\mu_{j-1}}{\longrightarrow} \ldots $$
$$ \ldots  \stackrel{\mu_{3}}{\longrightarrow} \kq_{\zeta}^{*}\otimes \bigwedge 
^{2j-1}\kq_{\zeta}^{*}(-2\gamma) \stackrel{\mu_{2}}{\longrightarrow} \kq_{\zeta}^{*}\otimes 
\bigwedge ^{2j}\kq_{\zeta}^{*}(-\gamma) \stackrel{\mu_{1}}{\longrightarrow} 
\kq_{\zeta}^{*}\otimes \bigwedge ^{2j+1}\kq_{\zeta}^{*}\stackrel{\mu_{0}}{\longrightarrow} 
\kq_{\zeta}^{*}\otimes \bigwedge ^{2j+1}\kn_{\zeta} \longrightarrow    0.$$

On consid\`ere  $\Sigma _{m}=ker (\mu_{m})$ pour $m=0,1, \ldots ,2j-1$. Les termes dans cette r\'esolution sont de la forme
$$\kq_{\zeta}^{*}\otimes\bigwedge ^{q}\kq_{\zeta}^{*}(-\alpha \gamma ),$$

o\`u $0\leq \alpha  \leq 2j ,\hspace{0.2 cm} 1\leq q\leq 2j+1 \leq n \text{ tels que } q+\alpha =2j+1 \text{ et } q \neq \alpha$. D'apr\`es la proposition \ref{3.5a}, on a 
$$h^{\alpha}(\kq_{\zeta}^{*}\otimes\bigwedge ^{2j+1-\alpha}\kq_{\zeta}^{*}(-\alpha \gamma ) )=0 \text{ pour tout } \alpha, \text{ si } \gamma > \sum_{i=0}^{n}\lambda_{i}.$$

Comme $h^{i}(\kq_{\zeta}^{*}(- \gamma(2j+1)))=0  \text{ pour tout } 0\leq i \leq 2n$, alors on obtient que $h^{1}(\Sigma_{0})=0$. En~utilisant la suite exacte suivante 
$$0\longrightarrow \Sigma_{0} \longrightarrow \kq_{\zeta}^{*}\otimes\bigwedge ^{2j+1}\kq_{\zeta}^{*}\stackrel{\mu_{0}}{\longrightarrow}\kq_{\zeta}^{*}\otimes\bigwedge ^{2j+1}\kn_{\zeta} \longrightarrow  0,$$

on obtient que 
$$h^{0}(\kq_{\zeta}^{*}\otimes \bigwedge^{2j+1}\kn_{\zeta})=0,$$

et
$$h^{0}(\kn_{\zeta}\otimes \bigwedge^{2j+1}\kn_{\zeta})=1 \hspace{0.2 cm}si\hspace{0.2 cm} 
 \gamma > \sum_{i=0}^{n}\lambda_{i}.$$

\end{proof}
   
\subsection{\bf Proposition}\label{3.6.1} 
{\em Soit $\kn_{\zeta}$ le fibr\'e de \mbox{0-corr\'elation} pond\'er\'e sur $\P^{2n+1}$ qui est d\'efini par la suite exacte \ref{8} pour $\zeta=1$. Soit $1\leq 2j+1\leq n$. Si on a $\gamma -n> \sum_{i=0}^{n}\lambda_{i}$, alors on obtient que 
$$h^{0}(\kn_{\zeta}^{*}(-(j+1))\otimes \bigwedge^{2j+1}\kn_{\zeta}^{*})=1.$$}

\begin{proof}
On fixe  $1\leq 2j+1\leq n$ et $\gamma -n > \sum_{i=0}^{n}\lambda_{i}$ et $\zeta=1$. D'apr\`es le lemme \ref{0.2.3} et la proposition \ref{2.3}, on obtient que  
$$h^{0}(\kn_{\zeta}^{*}(-(j+1))\otimes \bigwedge^{2j+1}\kn_{\zeta}^{*})\geq 1.$$ 

De la suite exacte suivante
$$0\longrightarrow \bigwedge^{2j+1}\kn_{\zeta}^{*}(-\gamma -(j+1))\longrightarrow \kq_{\zeta}^{*} (-(j+1))\otimes \bigwedge^{2j+1}\kn_{\zeta}^{*}\longrightarrow h^{0}(\kn_{\zeta}^{*}(-(j+1))\otimes \bigwedge^{2j+1}\kn_{\zeta}^{*})\longrightarrow 0$$

et suivant la proposition \ref{3.3a.1}, il suffit de d\'emontrer que 
$$h^{0}(\kq_{\zeta}^{*} (-(j+1))\otimes \bigwedge^{2j+1}\kn_{\zeta}^{*})=0$$

pour que l'on ait
$$h^{0}(\kn_{\zeta}^{*}(-(j+1))\otimes \bigwedge^{2j+1}\kn_{\zeta}^{*})=1.$$

De la r\'esolution \ref{6}, on obtient la r\'esolution suivante
$$0 \longrightarrow \kq_{\zeta}^{*}(- \gamma(2j+1)-(j+1))\longrightarrow \kq_{\zeta}^{*}\otimes \kq_{\zeta}^{*}(-\gamma(2j)-(j+1))\longrightarrow$$ 
$$\kq_{\zeta}^{*}\otimes \bigwedge ^{2}\kq_{\zeta}^{*}(-\gamma(2j-1)-(j+1)) \stackrel{\mu_{2j-1}}{\longrightarrow}
\kq_{\zeta}^{*}\otimes \bigwedge ^{3}\kq_{\zeta}^{*}(-\gamma(2j-2)-(j+1)) \stackrel{\mu_{2j-2}}{\longrightarrow}\ldots$$
$$\ldots \stackrel{\mu_{j+2}}{\longrightarrow} \kq_{\zeta}^{*}\otimes \bigwedge^{j} \kq_{\zeta}^{*}(- \gamma(j+1)-(j+1))\stackrel{\mu_{j+1}}{\longrightarrow} \kq_{\zeta}^{*}\otimes \bigwedge^{j+1}\kq_{\zeta}^{*}(-\gamma(j)-(j+1))\stackrel{\mu_{j}}{\longrightarrow} \ldots$$
 $$\ldots \stackrel{\mu_{1}}{\longrightarrow} \kq_{\zeta}^{*}\otimes \bigwedge ^{2j+1}\kq_{\zeta}^{*}(-(j+1))\stackrel{\mu_{0}}{\longrightarrow} \kq_{\zeta}^{*}(-(j+1))\otimes \bigwedge ^{2j+1}\kn_{\zeta}^{*} \longrightarrow    0$$

On consid\`ere  $\Sigma _{m}=Ker (\mu_{m})$ pour $m=0,1,\ldots ,2j-1$. Les termes dans cette r\'esolution sont de la forme
$$\kq_{\zeta}^{*}\otimes\bigwedge^{q}\kq_{\zeta}^{*}(-\alpha \gamma -(j+1)),$$

o\`u $0\leq \alpha  \leq 2j ,\hspace{0.2 cm} 1\leq q\leq 2j+1 \leq n $ tels que $q+\alpha =2j+1 \hspace{0.2 cm} q \neq \alpha$. D'apr\`es la proposition \ref{3.5a.1}; on a 
$$h^{\alpha}(\kq_{\zeta}^{*}\otimes\bigwedge ^{2j+1-\alpha}\kq_{\zeta}^{*}(-\alpha \gamma -(j+1)) )=0,\hspace{0.2 cm} pour\hspace{0.2 cm} tout\hspace{0.2 cm} \alpha$$

Comme $h^{i}(\kq_{\zeta}^{*}(- \gamma(2j+1))-(j+1))=0,\hspace{0.2 cm} pour\hspace{0.2 cm} tout\hspace{0.2 cm} 0\leq i \leq 2n$, alors on obtient que $h^{1}(\Sigma_{0})=0$. En utilisant de la suite exacte suivante 
$$0\longrightarrow \Sigma_{0} \longrightarrow \kq_{\zeta}^{*}\otimes\bigwedge ^{2j+1}\kq_{\zeta}^{*}(-(j+1)) \stackrel{\mu_{0}}{\longrightarrow}\kq_{\zeta}^{*}(-(j+1)) \otimes\bigwedge ^{2j+1}\kn_{\zeta}^{*} \longrightarrow  0$$

on obtient que 
$$h^{0}(\kq_{\zeta}^{*}(-(j+1)) \otimes\bigwedge ^{2j+1}\kn_{\zeta}^{*})=0.$$

Donc, on a  
$$ h^{0}(\kn_{\zeta}^{*}(-(j+1))\otimes \bigwedge^{2j+1}\kn_{\zeta}^{*})=1 \hspace{0.2 cm}si\hspace{0.2 cm} \gamma -n > \sum_{i=0}^{n}\lambda_{i}.$$

\end{proof}

\subsection{\bf Th\'eor\`eme}\label{3.7} 

{\em Soient $\kq_{\zeta;\gamma;\lambda_{n},\lambda_{n-1},\ldots ,\lambda_{0}}$ le fibr\'e de quotient pond\'er\'e sur $\P^{2n+1}$ des poids $\gamma;\lambda_{n},\lambda_{n-1}, \ldots ,\lambda_{0}$ et $\kn_{\zeta;\gamma;\lambda_{n},\lambda_{n-1}, \ldots ,\lambda_{0}}$ le fibr\'e de 0-corr\'elation pond\'er\'e sur $\P^{2n+1}$ des poids $\gamma;\lambda_{n},\lambda_{n-1}, \ldots ,\lambda_{0}$ qui sont d\'efinis par les suites exactes suivantes 
$$ 0\longrightarrow \ko_{\P^{2n+1}}(-\gamma)\stackrel{B}{\longrightarrow} \kh_{\zeta}\longrightarrow \kq_{\zeta;\gamma;\lambda_{n},\lambda_{n-1}, \ldots ,\lambda_{0}} \longrightarrow 0$$
$$0\longrightarrow \kn_{\zeta;\gamma;\lambda_{n},\lambda_{n-1}, \ldots ,\lambda_{0}}\longrightarrow  \kq_{\zeta;\gamma;\lambda_{n},\lambda_{n-1}, \ldots ,\lambda_{0}} \stackrel{A}{\longrightarrow} \ko_{\P^{2n+1}}(\gamma) \longrightarrow 0. $$

Pour $1\leq 2j+1\leq n$, les conditions suivantes sont \'equivalentes
 
I-  $\gamma-\zeta n > \sum_{i=0}^{n}\lambda_{i}$.

II- $\kn_{\zeta;\gamma;\lambda_{n},\lambda_{n-1}, \ldots ,\lambda_{0}}$ est stable.

III- $\kn_{\zeta;\gamma;\lambda_{n},\lambda_{n-1}, \ldots ,\lambda_{0}}$ est 
simple}.

\begin{proof}Pour $\zeta=0$.

$(I)\Longrightarrow (II)$. On fixe  $\gamma > \sum_{i=0}^{n}\lambda_{i}$. D'apr\`es la proposition \ref{2.3}, le fibr\'e $\kn_{\zeta;\gamma;\lambda_{n},\lambda_{n-1}, \ldots ,\lambda_{0}}$ est un fibr\'e symplectique. D'apr\`es les propositions \ref{3.3a} , \ref{3.6} et le th\'eor\`eme \ref{0.2.5}, on obtient que $\kn_{\zeta;\gamma;\lambda_{n},\lambda_{n-1}, \ldots ,\lambda_{0}}$ est stable.

$(II)\Longrightarrow (III)$. Evident.

$(III)\Longrightarrow (I)$. On suppose que $\gamma \leq 
\sum_{i=0}^{n}\lambda_{i}$. On choisit $\gamma_{0} \leq 
\sum_{i=0}^{1}\lambda_{n-i}$ un entier tel que $\gamma_{0} > \lambda_{n}\geq 
\lambda_{n-1}\geq \ldots \geq\lambda_{0}\geq 0$.  De la r\'esolution \ref{5}, 
on obtient la r\'esolution suivante
$$0 \longrightarrow \ko_{\P^{2n+1}}(-2n\gamma_{0})\longrightarrow 
\kh_{\zeta}(-\gamma_{0}(2n-1))\longrightarrow \bigwedge 
^{2}\kh_{\zeta}(-\gamma_{0}(2n-2))\stackrel{a_{2n-3}}{\longrightarrow}$$
$$\bigwedge 
^{3}\kh_{\zeta}(-\gamma_{0}(2n-3))\stackrel{a_{2n-4}}{\longrightarrow} \ldots $$
$$ \ldots  \stackrel{a_{3}}{\longrightarrow}\bigwedge ^{5}\kh_{\zeta}(-3\gamma_{0}) 
\stackrel{a_{2}}{\longrightarrow}\bigwedge ^{4}\kh_{\zeta}(-2\gamma_{0}) 
\stackrel{a_{1}}{\longrightarrow}\bigwedge 
^{3}\kh_{\zeta}(-\gamma_{0})\stackrel{a_{0}}{\longrightarrow}\bigwedge 
^{2}\kq_{\zeta;\gamma_{0};\lambda_{n},\lambda_{n-1}, \ldots ,\lambda_{0}}^{*} 
\longrightarrow   0.$$

On consid\`ere $A_{j}=ker (a_{j})$ pour $j=0,1, \ldots ,2n-3$. Du lemme \ref{3.1}, on obtient alors 
$$h^{k}(A_{0})=0,\hspace{0.2 cm} 0\leq k\leq 2.$$

De la suite exacte suivante
$$0\longrightarrow A_{0}\longrightarrow\bigwedge ^{3}\kh_{\zeta}(-\gamma_{0})\stackrel{a_{0}}{\longrightarrow}\bigwedge ^{2}\kq_{\zeta;\gamma_{0};\lambda_{n},\lambda_{n-1}, \ldots ,\lambda_{0}}^{*} \longrightarrow   0$$

on obtient la suite exacte des cohomologies 
$$0\longrightarrow H^{0}( A_{0})\longrightarrow H^{0}( \bigwedge ^{3}\kh_{\zeta}(-\gamma_{0}) )\stackrel{a_{0}}{\longrightarrow} H^{0}( \bigwedge ^{2}\kq_{\zeta;\gamma_{0};\lambda_{n},\lambda_{n-1}, \ldots ,\lambda_{0}}^{*} ) \longrightarrow   0.$$

D'apr\`es le lemme \ref{3.1} on obtient $ H^{0}( \bigwedge ^{3}\kh_{\zeta}(-\gamma_{0}) )\neq 0 $, car $\gamma_{0} \leq \sum_{i=0}^{1}\lambda_{n-i}$, ce qui nous donne
$$0\neq H^{0}( \bigwedge^{2}\kq_{\zeta;\gamma_{0};\lambda_{n},\lambda_{n-1}, \ldots ,\lambda_{0}}^{*})\neq \C.$$

De la r\'esolution \ref{6}, on obtient la r\'esolution suivante
$$0\longrightarrow  \ko_{\P^{2n+1}}(-2\gamma_{0})\longrightarrow  
\kq_{\zeta;\gamma_{0};\lambda_{n},\lambda_{n-1}, \ldots ,
\lambda_{0}}^{*}(-\gamma_{0} ) \longrightarrow  
\bigwedge^{2}\kq_{\zeta;\gamma_{0};\lambda_{n},\lambda_{n-1}, \ldots ,
\lambda_{0}}^{*}
\stackrel{b_{0}}{\longrightarrow}$$
$$\qquad\qquad\qquad  \bigwedge 
^{2}\kn_{\zeta;\gamma_{0};\lambda_{n},\lambda_{n-1}, \ldots ,\lambda_{0}}  
\longrightarrow   0$$

qui est \'equivalente aux deux suites exactes suivantes
$$0\longrightarrow B_{0} \longrightarrow  \bigwedge ^{2}\kq_{\zeta;\gamma_{0};\lambda_{n},\lambda_{n-1}, \ldots ,\lambda_{0}}^{*}\stackrel{b_{0}}{\longrightarrow}  \bigwedge ^{2}\kn_{\zeta;\gamma_{0};\lambda_{n},\lambda_{n-1}, \ldots ,\lambda_{0}}  \longrightarrow   0,$$
$$0\longrightarrow  \ko_{\P^{2n+1}}(-2\gamma_{0})\longrightarrow  \kq_{\zeta;\gamma_{0};\lambda_{n},\lambda_{n-1}, \ldots ,\lambda_{0}}^{*}(-\gamma_{0} ) \longrightarrow B_{0} \longrightarrow 0.$$

On obtient $h^{0}(B_{0})=0$, $h^{1}(B_{0})=1$. De la suite exacte de cohomologie des fibr\'es
$$ 0\longrightarrow H^{0}( \bigwedge ^{2}\kq_{\zeta;\gamma_{0};\lambda_{n},\lambda_{n-1}, \ldots ,\lambda_{0}}^{*})\longrightarrow H^{0}(\bigwedge ^{2}\kn_{\zeta;\gamma_{0};\lambda_{n},\lambda_{n-1}, \ldots ,\lambda_{0}}  )  \longrightarrow  H^{1}(B_{0})\longrightarrow 0,$$

on obtient que $ 0 \neq H^{0}(\bigwedge ^{2}\kn_{\zeta;\gamma_{0};\lambda_{n},\lambda_{n-1}, \ldots ,\lambda_{0}})\neq \C$. Mais, on a 
$$ H^{0}(\bigwedge ^{2}\kn_{\zeta;\gamma_{0};\lambda_{n},\lambda_{n-1}, \ldots ,\lambda_{0}})\subseteq H^{0}(\kn_{\zeta;\gamma_{0};\lambda_{n},\lambda_{n-1}, \ldots ,\lambda_{0}}\otimes \kn_{\zeta;\gamma_{0};\lambda_{n},\lambda_{n-1}, \ldots ,\lambda_{0}}).$$ 

Autrement dit, $\kn_{\zeta;\gamma_{0};\lambda_{n},\lambda_{n-1}, \ldots ,\lambda_{0}}$ n'est pas simple. Ce qui est une contradiction.
 
Pour $\zeta=1$, la d\'emonstration dans ce cas est tr\`es similaire \`a celle pour $\zeta=0$.
\end{proof}

\newpage

\section{ \large\bf D\'eformation miniversale de fibr\'e de 0-corr\'elation pond\'er\'e }
 
\vspace{1cm} 
Tout en fixant $\zeta$, nous allons d\'emontrer que les fibr\'es $\kq_{\zeta}$ et $\kn_{\zeta}$ sont invariants par rapport \`a~une d\'eformation miniversale. Nous allons aussi montrer que l'espace de Kuranishi du fibr\'e $\kn_{\zeta}$ est~lisse au point correspondant de $\kn_{\zeta}$. Pour $\zeta=0$, l'espace du module de fibr\'e stable $\kn_{\zeta}$ sur $\P^{2n+1}$ est s\'eparable topologiquement au point correspondant au fibr\'e $\kn_{\zeta}$.   

\subsection{\bf Th\'eor\`eme}\label{}
{\em (Hartshorne, \cite{ha10}). Si $F$ est un faisceau coh\'erent sur un sch\'ema projectif~$X$ sur un corps de base $K$ tel que $hdF\leq 1$, il existe un sch\'ema $Y=Spec(R)$ qui param\'etrise les~d\'eformations miniversales de $F$, o\`u $R$ est une $K$-alg\`ebre locale compl\`ete.}
 
\begin{proof}
Voir le th\'eor\`eme (19.1 \cite{ha10}).
\end{proof}

\subsection{\bf Th\'eor\`eme} 
{\em Soit $E$ un fibr\'e vectoriel sur l'espace projectif $\P^{n}$. Il existe $Kur(E)$, un espace de \mbox{Kuranishi} de $E$, qui est une base de la d\'eformation miniversale de $E$. Autrement dit, $Kur(E)$ param\'etrise toutes les d\'eformations miniversales de $E$. }

\begin{proof}
Voir l'article de M. Kuranishi \cite{kur}.
\end{proof}

Soit $e$ un point correspondant au fibr\'e $E$. Alors l'espace $Kur(E)$ est \'equip\'e d'une famille \mbox{universelle}. La fibre $(Kur(E),e)$, un espace topologique point\'e, est unique \`a un automorphisme pr\`es.

\subsection{\bf Lemme}\label{4.3} 
{\em Soient  $ \kq_{\zeta}^{'}$ et $ \kq_{\zeta}^{''}$ deux fibr\'es de quotient pond\'er\'e sur $\P^{2n+1}$ qui sont d\'efinis par les suites exactes suivantes
$$0\longrightarrow \ko_{\P^{2n+1}}(-\gamma-\zeta)\longrightarrow \kh_{\zeta}\stackrel{q_{1}}{\longrightarrow} \kq_{\zeta}^{'} \longrightarrow 0$$

$$0\longrightarrow \ko_{\P^{2n+1}}(-\gamma-\zeta)\longrightarrow \kh_{\zeta}\stackrel{q_{2}}{\longrightarrow} \kq_{\zeta}^{''} \longrightarrow 0$$

tels qu'il existe un morphisme $\psi : \kq_{\zeta}^{'} \longrightarrow \kq_{\zeta}^{''}$. Alors 
il existe un morphisme $\varphi : \kh_{\zeta}  \longrightarrow \kh_{\zeta} $ tel que 
$q_{2}\circ \varphi= \psi\circ q_{1}$}.

\begin{proof}
De la deuxi\`eme suite, on obtient la suite exacte suivante
$$0\longrightarrow Hom(\kh_{\zeta},\ko_{\P^{2n+1}}(-\gamma-\zeta)) \longrightarrow 
Hom(\kh_{\zeta},\kh_{\zeta})\stackrel{q_{2}\hspace{0.1 cm}\circ \hspace{0.1 cm} \bullet 
}{\longrightarrow} Hom(\kh_{\zeta},\kq_{\zeta}^{''}) $$ $$\longrightarrow 
Ext^{1}(\kh_{\zeta},\ko_{\P^{2n+1}}(-\gamma-\zeta)) \longrightarrow 0.$$
 
Comme on a  
$$Ext^{1}(\kh_{\zeta},\ko_{\P^{2n+1}}(-\gamma-\zeta))=H^{1}(\kh^{*}_{\zeta}(-\gamma-\zeta))=0,$$

alors on obtient que
$$Hom(\kh_{\zeta},\kh_{\zeta})\stackrel{q_{2}\hspace{0.1 cm}\circ \hspace{0.1 cm} \bullet}{\longrightarrow} Hom(\kh_{\zeta},\kq_{\zeta}^{''})\longrightarrow 0 $$

Comme on a le morphisme $\psi\circ q_{1}:\kh_{\zeta}\longrightarrow \kq_{\zeta}^{''} $ qui est 
dans $Hom(\kh_{\zeta},\kq_{\zeta}^{''})$, donc il existe un morphisme $\varphi: 
\kh_{\zeta}\longrightarrow\kh_{\zeta}$ tel que $q_{2}\circ \varphi= \psi\circ q_{1}$.
 
\end{proof}

\subsection{\bf Lemme}\label{4.4a} 
{\em Soient  $f,\hspace{0.2cm}f^{'}\in Hom(\ko_{\P^{2n+1}}(-\gamma-\zeta),\kh_{\zeta}) $ deux 
morphismes. Alors $f$ et $f'$ donnent le m\^eme \'el\'ement dans le sch\'ema 
$\QQuot_{\kh_{\zeta}/\P^{2n+1}}$ si et seulement s'il existe un isomorphisme \mbox{$g \in End( 
\ko_{\P^{2n+1}}(-\gamma-\zeta))$} tel que $f=f^{'}\circ g$}.

\begin{proof}
C'est la d\'efinition de sch\'ema $\QQuot_{\kh_{\zeta}/\P^{2n+1}}$.
  
\end{proof}

 Le th\'eor\`eme suivant est une g\'en\'eralisation du (\cite{an-ot93}, 
th\'eor\`eme 3.3) sur $\P^{2n+1}$.

\subsection{\bf Th\'eor\`eme }\label{4.5} 

{\em Soit $ \kq_{\zeta}{}_{,0}$ un fibr\'e de quotient pond\'er\'e sur $\P^{2n+1}$ qui est d\'efini par la suite exacte suivante
$$0\longrightarrow \ko_{\P^{2n+1}}(-\gamma-\zeta)\stackrel{x_{0}}{\longrightarrow} \kh_{\zeta}\longrightarrow \kq_{\zeta}{}_{,0} \longrightarrow 0$$

o\`u $x_{0}\in Hom(\ko_{\P^{2n+1}}(-\gamma-\zeta),\kh_{\zeta})$.
Alors chaque d\'eformation miniversale de fibr\'e $\kq_{\zeta}{}_{,0}$ est encore un fibr\'e de quotient pond\'er\'e sur $\P^{2n+1}$. L'espace de Kuranishi de $\kq_{\zeta}{}_{,0}$ est lisse au point correspondant de $\kq_{\zeta}{}_{,0}$.}

\begin{proof}
Soient $x_{0}\in Hom(\ko_{\P^{2n+1}}(-\gamma-\zeta),\kh_{\zeta})$ et $\kq_{\zeta}{}_{,0}=coker(x_{0})$ le fibr\'e vectoriel quotient pond\'er\'e correspondant au $x_{0}$. Soit $X$ un composant irr\'eductible de $\QQuot_{\kh_{\zeta}/\P^{2n+1}}$ tel que $x_{0} \in X$ et $y_{0}\in Kur(\kq_{\zeta}{}_{,0})$ correspondant au fibr\'e $\kq_{\zeta}{}_{,0}$. On a le morphisme de fibres suivant
$$\pi : (X,x_{0})\longrightarrow (Kur(\kq_{\zeta}{}_{,0}),y_{0})$$

D'apr\`es le th\'eor\`eme 3.12, page 137 \cite{qin}, on a 
$$dim_{y_{0}}(Kur(\kq_{\zeta}{}_{,0}))\geq  dim_{x_{0}}(X)- dim_{x_{0}}(\pi^{-1}(y_{0})).$$

D'apr\`es le th\'eor\`eme 2.2, page 126 \cite{qin}, on a $h^{1}(\EEnd(\kq_{\zeta}{}_{,0})))\geq dim_{y_{0}}(Kur(\kq_{\zeta}{}_{,0}))$. Soit 
$$Z=\lbrace x_{1}\in X | \hspace{0.2 cm}\kq_{\zeta}{}_{,1}\simeq \kq_{\zeta}{}_{,0} \hspace{0.2cm} \text{o\`u $\kq_{\zeta}{}_{,1}$ est le fibr\'e correspondant \`a $x_{1}$ } \rbrace.$$ 

On obtient que
$$(\pi^{-1}(y_{0}),y_{0})\subseteq (Z,x_{0})\hspace{0.2cm} et\hspace{0.2cm} dim_{x_{0}}((\pi^{-1}(y_{0}),y_{0}))\leq dim_{x_{0}}((Z,x_{0})).$$

Donc on a 
$$dim_{y_{0}}(Kur(\kq_{\zeta}{}_{,0}))\geq  dim_{x_{0}}(X)-dim_{x_{0}}((Z,x_{0}))  .$$

On a aussi $ dim_{x_{0}}(X)=h^{0}(\kh_{\zeta}(\gamma+\zeta))-1=h^{0}(\kq_{\zeta}{}_{,0}(\gamma+\zeta))$. Soit 
$$\Sigma =\lbrace \psi \in End(\kh_{\zeta})| \hspace{0.2 cm}\psi.x_{0}=x_{0}\rbrace.$$ 

D'apr\`es les lemmes \ref{4.3} et \ref{4.4a}, on obtient que 
$$dim_{x_{0}}(Z)= h^{0}(\EEnd(\kh_{\zeta}))-dim_{x_{0}}(\Sigma)-1.$$

De la suite exacte de fibr\'es vectoriels
$$0\longrightarrow \kq_{\zeta}^{*}{}_{,0}\otimes \kh_{\zeta} \longrightarrow \EEnd(\kh_{\zeta})\stackrel{}{\longrightarrow} \kh_{\zeta}(\gamma+\zeta) \longrightarrow 0,$$

on obtient la suite exacte suivante de groupes cohomologiques
$$0\longrightarrow H^{0}(\kq_{\zeta}{}^{*}_{,0}\otimes \kh_{\zeta}) \longrightarrow End(\kh_{\zeta})\stackrel{\bullet\circ x_{0}}{\longrightarrow} H^{0}(\kh_{\zeta}(\gamma+\zeta)) ,$$

qui nous donne
$$dim_{x_{0}}(\Sigma)=h^{0}(\kq_{\zeta}{}^{*}_{,0}\otimes \kh_{\zeta}).$$ 

De la suite exacte pr\'ec\'edente de fibr\'es vectoriels, on obtient que 
$$h^{1}(\kq_{\zeta}{}^{*}_{,0}\otimes \kh_{\zeta})=h^{0}(\kh_{\zeta}(\gamma+\zeta))- h^{0}( \EEnd(\kh_{\zeta}))+ h^{0}(\kq_{\zeta}{}^{*}_{,0}\otimes \kh_{\zeta}).  $$

Donc on a
$$dim_{y_{0}}(Kur(\kq_{\zeta}{}_{,0}))\geq  h^{0}(\kh_{\zeta}(\gamma+\zeta))- h^{0}( \EEnd(\kh_{\zeta}))+ h^{0}(\kq_{\zeta}{}^{*}_{,0}\otimes \kh_{\zeta}) = h^{1}(\kq_{\zeta}{}^{*}_{,0}\otimes \kh_{\zeta}).$$

La suite exacte suivante des fibr\'es vectoriels
$$0\longrightarrow \kq_{\zeta}{}^{*}_{,0}(-\gamma-\zeta)\longrightarrow \kq_{\zeta}{}^{*}_{,0}\otimes\kh_{\zeta}\longrightarrow \EEnd(\kq_{\zeta}{}_{,0}) \longrightarrow 0$$

nous donne que $h^{2}(\EEnd(\kq_{\zeta}{}_{,0}))=0$ et $H^{1}(\kq_{\zeta}{}^{*}_{,0}\otimes\kh_{\zeta})\longrightarrow H^{1}(\EEnd(\kq_{\zeta}{}_{,0}))\longrightarrow 0$. Donc on a \mbox{$h^{1}(\kq_{\zeta}{}^{*}_{,0}\otimes\kh_{\zeta})\geq h^{1}(\EEnd(\kq_{\zeta}{}_{,0}))$}. Ensuite on obtient que
$$dim_{y_{0}}(Kur(\kq_{\zeta}{}_{,0}))\geq h^{1}(\EEnd(\kq_{\zeta}{}_{,0})).$$

Autrement dit, $dim_{y_{0}}(Kur(\kq_{\zeta}{}_{,0}))=h^{1}(\EEnd(\kq_{\zeta}{}_{,0}))$ et $Kur(\kq_{\zeta}{}_{,0})$ est lisse en $y_{0}$. De plus on a
$$dim_{y_{0}}(Kur(\kq_{\zeta}{}_{,0}))=  dim_{x_{0}}(X)- dim_{x_{0}}(\pi^{-1}(y_{0})).$$

D'apr\`es le th\'eor\`eme de la semicontinuit\'e de fibres (12.8, page 288 \cite{ha}), on obtient que\\ $dim_{y_{0}}(Im(\pi))= dim_{y_{0}}(Kur(\kq_{\zeta}{}_{,0}))$ et que $\pi $ est surjectif. Cela implique que $\kq_{\zeta}{}_{,0}$ est invariant par rapport \`a une d\'eformation miniversale.
 
\end{proof}

\subsection{\bf Lemme}\label{4.6} 
{\em Soit $\kq_{\zeta}$ un fibr\'e de quotient pond\'er\'e sur $\P^{2n+1}$. Soient $\kn_{\zeta}^{'}{}^{*}$ et $ \kn_{\zeta}^{''}{}^{*}$ des fibr\'es de \mbox{0-corr\'elation} pond\'er\'e sur $\P^{2n+1}$ qui sont d\'efinis par les suites exactes
$$0\longrightarrow \ko_{\P^{2n+1}}(-\gamma) \longrightarrow \kq_{\zeta}^{*}  \stackrel{p_{1}}{\longrightarrow} \kn_{\zeta}^{'}{}^{*} \longrightarrow 0,$$

et
$$0\longrightarrow \ko_{\P^{2n+1}}(-\gamma) \longrightarrow \kq_{\zeta}^{*}  \stackrel{p_{2}}{\longrightarrow} \kn_{\zeta}^{''}{}^{*} \longrightarrow 0.$$

tels qu'il existe un morphisme $\psi :\kn_{\zeta}^{'}{}^{*} \longrightarrow \kn_{\zeta}^{''}{}^{*} 
$. Alors il existe un morphisme $\varphi :\kq_{\zeta}^{*} \longrightarrow 
\kq_{\zeta}^{*}  $ tel que $p_{2}\circ \varphi= \psi\circ p_{1}$.}

\begin{proof}
 De la suite exacte suivante
$$0\longrightarrow \ko_{\P^{2n+1}}(-\gamma) \longrightarrow \kq_{\zeta}^{*}  \stackrel{p_{2}}{\longrightarrow} \kn_{\zeta}^{''}{}^{*} \longrightarrow 0,$$

on obtient la suite exacte suivante de groupes cohomologiques
$$0\longrightarrow Hom(\kq_{\zeta}^{*},\ko_{\P^{2n+1}}(-\gamma)) \longrightarrow 
Hom(\kq_{\zeta}^{*},\kq_{\zeta}^{*})\stackrel{p_{2}\hspace{0.1 cm}\circ \hspace{0.1 
cm} \bullet }{\longrightarrow} Hom(\kq_{\zeta}^{*},\kn_{\zeta}^{''}{}^{*}) $$
$$\longrightarrow Ext^{1}(\kq_{\zeta}^{*},\ko_{\P^{2n+1}}(-\gamma)) 
\longrightarrow 0.$$
 
Comme on a  
$$Ext^{1}(\kq_{\zeta}^{*},\ko_{\P^{2n+1}}(-\gamma))=H^{1}(\kq_{\zeta}(-\gamma))=0,$$

alors on obtient que 
$$Hom(\kq_{\zeta}^{*},\kq_{\zeta}^{*})\stackrel{p_{2}\hspace{0.1 cm}\circ \hspace{0.1cm} \bullet }{\longrightarrow} Hom(\kq_{\zeta}^{*},\kn_{\zeta}^{''}{}^{*})\longrightarrow 0.$$

Comme on a le morphisme $\psi\circ p_{1}:\kq_{\zeta}^{*}\longrightarrow \kn_{\zeta}^{''}{}^{*} $ qui est dans $Hom(\kq_{\zeta}^{*},\kn_{\zeta}^{''}{}^{*})$, donc il existe un 
morphisme $\varphi: \kq_{\zeta}^{*}\longrightarrow \kq_{\zeta}^{*}$ tel que 
$p_{2}\circ \varphi= \psi\circ p_{1}$.
 
\end{proof}

\subsection{\bf Lemme}\label{4.7}
{\em Soit $\kq_{\zeta}$ un fibr\'e de quotient pond\'er\'e sur $\P^{2n+1}$. Soient $f, f^{'}$ deux morphismes dans $Hom(\ko_{\P^{2n+1}}(-\gamma),\kq_{\zeta}^{*}) $. Alors $f$ et $f^{'}$ donnent le m\^eme \'el\'ement dans $\QQuot_{\kq_{\zeta}^{*}/\P^{2n+1}}$ si~et seulement s'il existe un isomorphisme $g \in End( \ko_{\P^{2n+1}}(-\gamma))$ tel que $f=f^{'}\circ g$}.

\begin{proof}
C'est la d\'efinition de sch\'ema $\QQuot_{\kq_{\zeta}^{*}/\P^{2n+1}}$.
 
\end{proof}

Le th\'eor\`eme suivant est une g\'en\'eralisation du (\cite{an-ot93}, 
th\'eor\`eme 4.4) sur $\P^{2n+1}$.
 
\subsection{\bf Th\'eor\`eme} \label{4.8} 

{\em Soient $\kq_{\zeta}$ le fibr\'e de quotient pond\'er\'e sur $\P^{2n+1}$ et $\kn_{\zeta}{}_{,0}$ le fibr\'e de \mbox{0-corr\'elation} pond\'er\'e sur $\P^{2n+1}$ qui sont d\'efinis par les suites exactes 
$$0\longrightarrow \ko_{\P^{2n+1}}(-\gamma) \stackrel{g_{0}}{\longrightarrow} \kq_{\zeta}^{*}  \stackrel{}{\longrightarrow} \kn_{\zeta}^{*}{}_{,0} \longrightarrow 0$$

et
$$0\longrightarrow \ko_{\P^{2n+1}}(-\gamma-\zeta)\stackrel{}{\longrightarrow} \kh_{\zeta}\longrightarrow \kq_{\zeta} \longrightarrow 0,$$

o\`u $g_{0}\in Hom(\ko_{\P^{2n+1}}(-\gamma),\kq_{\zeta}^{*})$. Alors chaque d\'eformation miniversale de fibr\'e $\kn_{\zeta}{}_{,0}$ est encore un fibr\'e de 0-corr\'elation pond\'er\'e sur $\P^{2n+1}$. L'espace de Kuranishi de $\kn_{\zeta}{}_{,0}$ est lisse au point correspondant de $\kn_{\zeta}{}_{,0}$.}

\begin{proof}
Soient $g_{0}\in Hom(\ko_{\P^{2n+1}}(-\gamma),\kq_{\zeta}^{*}) $ et $\kn_{\zeta}^{*}{}_{,0}=coker(g_{0})$ un fibr\'e vectoriel \mbox{quotient} de $\kq_{\zeta}^{*}$ correspondant au $g_{0}$. Soit $Y \subseteq \QQuot_{\kq_{\zeta}^{*}/\P^{2n+1}}$ un composant irr\'eductible de $\QQuot_{\kq_{\zeta}^{*}/\P^{2n+1}}$ tel que  $g_{0} \in Y$. Soient $x \in Kur(\kq_{\zeta})$ correspondant au fibr\'e $\kq_{\zeta}$ et $z_{0} \in Kur(\kn_{\zeta}{}_{,0})$ correspondant au fibr\'e $\kn_{\zeta}{}_{,0}$.

\xmat{  (Y,g_{0}) \ar[rr]^{\Psi}\ar[dd]_{\Phi} &&(Kur\kq_{\zeta},x)\\
\\ 
(Kur\kn_{\zeta}{}_{,0},z_{0})&&\\}

D'apr\`es le th\'eor\`eme \ref{4.5}, pour le morphisme $\Psi$, on a
$$  dim_{g_{0}}(Y)=dim_{x}(Kur(\kq_{\zeta}))+ dim_{g_{0}}(\Psi^{-1}(x)),$$

et $dim_{x}(Kur(\kq_{\zeta}))= h^{1}(\EEnd(\kq_{\zeta}))$. La dimension de la fibre du morphisme $\Psi$ est \'egale \`a $h^{0}(\kq_{\zeta}^{*}(\gamma))-h^{0}(\ko_{\P^{2n+1}})$, donc on obtient 
$$  dim_{g_{0}}(Y)=h^{1}(\EEnd(\kq_{\zeta}))+ h^{0}(\kq_{\zeta}^{*}(\gamma))-1.$$

En faisant la m\^eme chose comme dans le th\'eor\`eme \ref{4.5} on obtient, pour le morphisme $\Phi$,
$$h^{1}(\EEnd(\kn_{\zeta}{}_{,0}))\geq dim_{z_{0}}(Kur(\kn_{\zeta}{}_{,0}))\geq  dim_{g_{0}}(Y)- dim_{g_{0}}(\Phi^{-1}(z_{0})).$$

Soit $Z=\lbrace g_{1}\in Y | \hspace{0.2 cm}\kn_{\zeta}^{*}{}_{,1}\simeq \kn_{\zeta}^{*}{}_{,0} \hspace{0.2cm} \text{o\`u $\kn_{\zeta}^{*}{}_{,1}$ est le fibr\'e correspondant \`a $g_{1}$ } \rbrace$. On obtient que
$$(\Phi^{-1}(z_{0}),g_{0})\subseteq (Z,g_{0})\hspace{0.2cm} et\hspace{0.2cm} dim_{g_{0}}((\Phi^{-1}(z_{0}),g_{0}))\leq dim_{g_{0}}((Z,g_{0})).$$

Donc on a 
$$dim_{z_{0}}(Kur(\kn_{\zeta}{}_{,0}))\geq  dim_{g_{0}}(Y)-dim_{g_{0}}((Z,g_{0})).$$

Soit $\Sigma =\lbrace \sigma \in End(\kq_{\zeta})| \hspace{0.2 cm}\sigma .g_{0}=g_{0}\rbrace$. D'apr\`es les lemmes \ref{4.6} et \ref{4.7}, on obtient alors que
$$dim_{g_{0}}(Z)= h^{0}(\EEnd(\kq_{\zeta}))-dim_{g_{0}}(\Sigma)-1.$$

En consid\'erant la suite exacte suivante de fibr\'es vectoriels
$$0\longrightarrow \kq_{\zeta}^{*}\otimes \kn_{\zeta}{}_{,0}\longrightarrow \EEnd(\kq_{\zeta})\stackrel{}{\longrightarrow} \kq_{\zeta}^{*}(\gamma) \longrightarrow 0,$$

on obtient la suite exacte suivante de groupes cohomologiques 
$$0\longrightarrow H^{0}(\kq_{\zeta}^{*}\otimes \kn_{\zeta}{}_{,0}) \longrightarrow End(\kq_{\zeta})\stackrel{\bullet\circ g_{0}}{\longrightarrow} H^{0}(\kq_{\zeta}^{*}(\gamma)) ,$$

qui nous donne $dim_{g_{0}}(\Sigma)=h^{0}(\kq_{\zeta}^{*}\otimes \kn_{\zeta}{}_{,0})$. Donc on obtient que 
$$dim_{z_{0}}(Kur(\kn_{\zeta}{}_{,0}))\geq h^{1}(\EEnd(\kq_{\zeta}))+h^{0}(\kq_{\zeta}^{*}(\gamma))- h^{0}(\EEnd(\kq_{\zeta}))+h^{0}(\kq_{\zeta}^{*}\otimes \kn_{\zeta}{}_{,0}).$$

De la suite exacte pr\'ec\'edente de fibr\'es vectoriels, on obtient que 
$$dim_{z_{0}}(Kur(\kn_{\zeta}{}_{,0}))\geq h^{1}(\EEnd(\kq_{\zeta}))+h^{0}(\kq_{\zeta}^{*}(\gamma))- h^{0}(\EEnd(\kq_{\zeta}))+h^{0}(\kq_{\zeta}^{*}\otimes \kn_{\zeta}{}_{,0})=$$
$$=h^{1}(\kq_{\zeta}^{*}\otimes \kn_{\zeta}{}_{,0})-h^{2}(\kq_{\zeta}^{*}\otimes \kn_{\zeta}{}_{,0})+h^{1}(\kq_{\zeta}^{*}(\gamma))+ h^{2}(\EEnd(\kq_{\zeta})) .$$

et 
$$\ldots \longrightarrow  H^{1}(\kq_{\zeta}^{*}(\gamma)) \longrightarrow H^{2}(\kq_{\zeta}^{*}\otimes \kn_{\zeta}{}_{,0}) \longrightarrow H^{2}(\EEnd(\kq_{\zeta}))\longrightarrow \ldots $$

Donc on obtient que $h^{2}(\kq_{\zeta}^{*}\otimes \kn_{\zeta}{}_{,0})\leq h^{2}(\EEnd(\kq_{\zeta}))+ h^{1}(\kq_{\zeta}^{*}(\gamma))$ et que
$$dim_{z_{0}}(Kur(\kn_{\zeta}{}_{,0}))\geq h^{1}(\kq_{\zeta}^{*}\otimes \kn_{\zeta}{}_{,0}).$$

En consid\'erant la suite exacte suivante de fibr\'es vectoriels
$$0\longrightarrow  \kn_{\zeta}{}_{,0}(-\gamma)  \longrightarrow  \kq_{\zeta}^{*}\otimes \kn_{\zeta}{}_{,0}  \longrightarrow \EEnd(\kn_{\zeta}{}_{,0}) \longrightarrow 0,$$

on obtient que 
$$\ldots \longrightarrow  H^{1}(\kq_{\zeta}^{*}\otimes \kn_{\zeta}{}_{,0})  \longrightarrow   H^{1}(\EEnd(\kn_{\zeta}{}_{,0}))  \longrightarrow H^{2}(\kn_{\zeta}{}_{,0}(-\gamma)) \longrightarrow \ldots $$

Donc on a
$$h^{1}(\EEnd(\kn_{\zeta}{}_{,0}))\leq h^{1}(\kq_{\zeta}^{*}\otimes \kn_{\zeta}{}_{,0})+ h^{2}(\kn_{\zeta}{}_{,0}(-\gamma))=h^{1}(\kq_{\zeta}^{*}\otimes \kn_{\zeta}{}_{,0}) $$

et
$$dim_{z_{0}}(Kur(\kn_{\zeta}{}_{,0}))\geq h^{1}(\EEnd(\kn_{\zeta}{}_{,0})).$$

Alors on a $dim_{z_{0}}(Kur(\kn_{\zeta}{}_{,0}))= h^{1}(\EEnd(\kn_{\zeta}{}_{,0}))$ et $Kur(\kn_{\zeta}{}_{,0})$ est lisse en $z_{0}$. De plus on a
$$dim_{z_{0}}(Kur(\kn_{\zeta}{}_{,0}))= dim_{g_{0}}(Y)- dim_{g_{0}}(\Phi^{-1}(z_{0})).$$

D'apr\`es le th\'eor\`eme de la semicontinuit\'e de fibres (12.8, page 288 \cite{ha}), on obtient que\\ $dim_{z_{0}}(Im(\Phi))= dim_{z_{0}}(Kur(\kn_{\zeta}{}_{,0}))$ et que $\Phi $ est surjectif. Cela implique que $\kn_{\zeta}{}_{,0}$ est invariant par rapport \`a une d\'eformation miniversale.
 
\end{proof}

\subsection{\bf Th\'eor\`eme} {\em (\cite{ko-ok}, th\'eor\`eme 6.4). \label{4.9} 
 Soient $E$ et $ E^{'}$ deux fibr\'es vectoriels simples non-isomorphiques
sur $\P^{m}$. Si les points associ\'es aux fibr\'es $E$ et $ E^{'}$ sont 
 non-s\'eparables (topologiquement) dans l'espace du module de fibr\'es 
 simples, alors il existe deux morphismes non-triviaux 
$$\varphi:E\longrightarrow E^{'},\hspace{0.2cm} \psi:E^{'}\longrightarrow E$$

tels que $\varphi\circ\psi=\psi\circ\varphi=0$.}

 \begin{proof}
 voir \cite{ko-ok}.
\end{proof}

\subsection{\bf Proposition} {\em (\cite{ok-sc-sp}, Lemme 1.2.8). \label{4.10} 
Soient $E$ et $ E^{'}$ deux fibr\'es vectoriels semi-stables tels que 
$rg(E^{'})=rg(E)$ et $c_{1}(E^{'})=c_{1}(E)$ sur $\P^{m}$. Soit 
$\varphi:E\longrightarrow E^{'}$ un morphisme non-trivial. Si au moins un des 
deux fibr\'es est stable, alors $ \varphi$ est un isomorphisme.}

 \begin{proof}
 voir \cite{ok-sc-sp}.
\end{proof}

\subsection{\bf Proposition}\label{4.11} 
{\em Soient $\kq_{\zeta}, \kq_{\zeta}^{'}$ des fibr\'es de quotient pond\'er\'e sur $\P^{2n+1}$ et $\kn_{\zeta},\kn_{\zeta}^{'}$ des fibr\'es de 0-corr\'elation pond\'er\'e sur $\P^{2n+1}$ qui sont d\'efinis par les suites exactes, pour $\zeta=0$, 

$$0\longrightarrow \ko_{\P^{2n+1}}(-\gamma)\longrightarrow \kh_{\zeta}\longrightarrow \kq_{\zeta} \longrightarrow 0$$
$$0\longrightarrow \ko_{\P^{2n+1}}(-\gamma)\longrightarrow \kh_{\zeta}\longrightarrow \kq_{\zeta}^{'} \longrightarrow 0,$$ 

et 
$$0\longrightarrow \ko_{\P^{2n+1}}(-\gamma)\longrightarrow \kq_{\zeta}^{*}\stackrel{q}{\longrightarrow}  \kn_{\zeta}\longrightarrow 0$$
$$0\longrightarrow \ko_{\P^{2n+1}}(-\gamma)\longrightarrow 
\kq_{\zeta}^{'}{}^{*}\stackrel{q_{'}}{\longrightarrow}  \kn_{\zeta}^{'}\longrightarrow 0.$$

Soit $\gamma > (2n+1)\lambda_{n}$. Alors les points associ\'es aux  
$\kn_{\zeta},\kn_{\zeta}^{'}$ sont s\'eparables dans l'espace du module 
$\km_{\P^{2n+1}}$.}
 
\begin{proof}
   On fixe $\gamma > (2n+1)\lambda_{n} \geq \sum_{i=0}^{n}\lambda_{i}$.
Supposons que les points associ\'es aux $\kn_{\zeta},\kn_{\zeta}^{'}$ sont non-s\'eparables dans l'espace de module $\km_{\P^{2n+1}}$. D'apr\`es la proposition \ref{4.9}, il existe deux morphismes non-triviaux 
$$\varphi:\kn_{\zeta}\longrightarrow \kn_{\zeta}^{'},\hspace{0.2cm} \psi:\kn_{\zeta}^{'}\longrightarrow \kn_{\zeta}$$

tels que $\varphi \circ\psi=\psi\circ\varphi=0$. En utilisant la suite exacte suivante
$$0\longrightarrow \ko_{\P^{2n+1}}(-\gamma)\longrightarrow 
\kq_{\zeta}^{'}{}^{*}\stackrel{q_{'}}{\longrightarrow}  \kn_{\zeta}^{'}\longrightarrow 0,$$

on obtient la suite exacte suivante de groupes cohomologiques 
$$0\longrightarrow Hom(\kq_{\zeta}^{*},\ko_{\P^{2n+1}}(-\gamma))\longrightarrow 
Hom(\kq_{\zeta}^{*}, \kq_{\zeta}^{'}{}^{*})\longrightarrow Hom(\kq_{\zeta}^{*}, 
\kn_{\zeta}^{'})$$
$$\longrightarrow Ext^{1}(\kq_{\zeta}^{*},\ko_{\P^{2n+1}}(-\gamma)).$$

Comme on a
$$Ext^{1}(\kq_{\zeta}^{*},\ko_{\P^{2n+1}}(-\gamma))=H^{1}(\kq_{\zeta}(-\gamma))=0,$$

alors on obtient
$$Hom(\kq_{\zeta}^{*}, \kq_{\zeta}^{'}{}^{*})\longrightarrow Hom(\kq_{\zeta}^{*}, 
\kn_{\zeta}^{'})\longrightarrow 0.$$  

Donc, pour le morphisme
$$\varphi \circ q: \kq_{\zeta}^{*}\stackrel{q}{\longrightarrow}  \kn_{\zeta} \stackrel{\varphi}{\longrightarrow} \kn_{\zeta}^{'} $$

il existe $0\neq \rho\in Hom(\kq_{\zeta}^{*}, \kq_{\zeta}^{'}{}^{*})$ tel que 
$q^{'}\circ \rho =\varphi \circ q$. D'apr\`es la proposition \ref{2.4} 
$\kq_{\zeta}^{'}{}^{*},\kq_{\zeta}^{*}$ sont stables. D'apr\`es la proposition \ref{4.10}, 
on obtient que le morphisme $\rho$ est un isomorphisme. Donc on a le 
diagramme commutatif

\xmat{    \kq_{\zeta}^{*} \ar[r]^{q}\ar[d]_{\rho}^{\wr}  & \kn_{\zeta} 
\ar[r]\ar[d]^{\varphi}& 0\\
   \kq_{\zeta}^{'}{}^{*} \ar[r]_{q^{'}} &  \kn_{\zeta}^{'} \ar[r]& 0\\
}

En r\'ep\'etant les m\^emes proc\'edures pr\'ec\'edentes pour le morphisme $\psi:\kn_{\zeta}^{'}\longrightarrow \kn_{\zeta}$, tout en utilisant la suite exacte

$$0\longrightarrow \ko_{\P^{2n+1}}(-\gamma)\longrightarrow \kq_{\zeta}^{*}\stackrel{q}{\longrightarrow}  \kn_{\zeta}\longrightarrow 0,$$

on obtient alors un isomorphisme $0\neq \rho^{'}\in Hom( {\kq_{\zeta}^{*}}^{'} 
,\kq_{\zeta}^{*})$ tel que le carr\'e suivant est commutatif

\xmat{    \kq_{\zeta}^{'}{}^{*} \ar[r]^{q^{'}}\ar[d]_{\rho^{'}}^{\wr}  & \kn_{\zeta}^{'} 
\ar[r]\ar[d]^{\psi}& 0\\
\kq_{\zeta}^{*} \ar[r]_{q} &  \kn_{\zeta}  \ar[r]& 0.\\  }

Donc on a le diagramme suivant

\xmat{ 0 \ar[r] & \ko_{\P^{2n+1}}(-\gamma) \ar[r] & 
\kq_{\zeta}^{*}\ar[r]^{q}\ar[d]_{\rho}^{\wr}  & \kn_{\zeta} \ar[r]\ar[d]^{\varphi}& 0\\
0 \ar[r] & \ko_{\P^{2n+1}}(-\gamma) \ar[r] & \kq_{\zeta}^{'*} \ar[r]_{q^{'}}\ar[d]_{\rho^{'}}^{\wr} &  \kn_{\zeta}^{'} \ar[r]\ar[d]^{\psi}& 0\\
0 \ar[r] & \ko_{\P^{2n+1}}(-\gamma) \ar[r] &\kq_{\zeta}^{*} \ar[r]_{q} &  \kn_{\zeta}  
\ar[r]& 0\\ }

tel que $q\circ\rho^{'} \circ \rho =\psi \circ \varphi \circ q=0$. Donc on a le morphisme suivant qui est d\'efini, pour tout $x\in \P^{2n+1} $, par

$$\hspace{0.7cm}g_{x}:\kq_{\zeta}{}^{*}_{x} \longrightarrow (\ko_{\P^{2n+1}}(-\gamma))_{x}\longrightarrow 0 $$
$$z  \longmapsto \rho^{'}\circ \rho(z).$$

Alors on a le morphisme $  0\longrightarrow \ko_{\P^{2n+1}} 
\stackrel{g^{*}}{\longrightarrow} \kq_{\zeta}(-\gamma) $, c'est-\`a-dire on a une section 
de fibr\'e $\kq_{\zeta}(-\gamma)$. Ce qui est une contradiction au fait que 
$H^{0}(\kq_{\zeta}(-\gamma))=0$.

\end{proof}


\newpage

\begin{thebibliography}{99}

\bibitem{an-ot94}Ancona, V. Ottaviani, G. {\em Stability of special instanton bundles on $\P^{2n+1}$}. Trans. of the Amer. Math. Soc.341, 2 (1994), 677-693.

\bibitem{an-ot93}Ancona, V. Ottaviani, G. {\em  3- bundles on $\P^{5}$}. (1993).

\bibitem{an-ot95}Ancona, V. Ottaviani, G. {\em On moduli of instanton bundles on $\P^{2n+1}$}. Pacific J. of Math. 171, 2 (1995), 343-351.

\bibitem{ati}  Atiyah, M.F. {\em Geometry on YangMills fields}, Scuola Normale Superiore Pisa, Pisa, 1979,
99 pp.

\bibitem{atwa}  Atiyah, M.F. R.S. Ward. {\em Instantons and algebraic geometry}, Comm. Math. Phys. 55
(1977), 117-124.

\bibitem{atdrhima} Atiyah, M. Drinfeld, V. Hitchin, N. and Manin, Y. {\em Construction of instantons}, Phys. Lett., 65A(1978), \mbox{pp. 185-187}.

\bibitem{bo-sp}Bohnhorst, G. Spindler, H. {\em The stability of certain vector bundles on $\P^{n}$}. Complex algebraic varieties (Bayreuth, 1990), 39-50, Lect. Notes in Math. , 1507, Springer, Berlin, 1992.

\bibitem{br} Br\^inz\u anescu, V. {\em Holomorphic Vector Bundles over
Compact Complex Surfaces.} Lect. Notes in Math. 1624. Springer-Verlag, Berlin (1996).

\bibitem{ca} Cascini, P. {\em Weighted Tango bundles on $\P^{n}$ and their moduli spaces}. Forum Math. 13(2001), 251-260.

\bibitem{dio} Dionisi, C. {\em Symplectic Small Deformations of Special Instanton Bundle on $\P^{2n+1}$}. Ann. di Mat. pura et appl. 175 (1998), 285-293.

\bibitem{dove} Douady, A. Verdier, J.-L. editors. {\em Les \'equations de Yang-Mills. S\'eminaire
E.N.S}. 1977-1978, Ast\'erisque 71-72. Paris: Soc. Math. France, 1980.

\bibitem{ei}Ein, L. {\em generalized null correlation bundles}. Nagoya Math. J.Vol. HI (1988), 13-24.

\bibitem{fu} Fulton, W. {\em Intersection Theory.} Springer-Verlag, Berlin (1998).

\bibitem{ha} Hartshorne, R. {\em Algebraic Geometry.} Gradua.Texte in Math. 52. 
Springer-Verlag, Berlin (1977).

\bibitem{ha10} Hartshorne, R. {\em Deformation theory.} Gradua.Texte in Math. 257. 
Springer New York (2010).

\bibitem{hu} Husemoller, D. {\em Fibre Bundles}. Third edition. Grad. Texts in
Math. 20. Springer-Verlag, New York (1994).

\bibitem{hule} Huybrecht, D. lehn, M. {\em The geometry of moduli space of scheaves}. Seco. edition. Cambr.Univ.Press. 2010. 

\bibitem{ja-mi}Jardim, M. Mir\`o-Roig, R. M. {\em  On the semistability of instanton sheaves over certain projective varieties}. Comm. in Algebra 36 (2008), 288-298.

\bibitem{ko-ok} Kosarew, S. Okonek, C. {\em Global Moduli Spaces and Simple Holomorphic Bundles}. Publ. RIMS, Kyoto Univ.25 (1989), 1-19.

\bibitem{kur} Kuranishi, M. {\em New proof for the existence of locally complete 
families of complex structures}. Proceedings of the Conference on Complex Analysis 1965, pp 142-154.

\bibitem{lep} Le Potier, J. {\em Lectures on vector bundles}. Cambridge Studies
in Adv. Math. 54. Cambridge University Press (1997).

\bibitem{mi-na-pe} Migliore, C. Nagel, U. Peterson, C. {\em Buchsbaum-Rim Sheaves and Their Multiple Sections}. Journal of \mbox{Algebra} 219, 378- 420. (1999).

\bibitem{ok-sc-sp}Okonek, C. Schneider, M. Spindler, H. {\em Vector Bundles on Complex Projective Spaces With an Appendix by S. I. Gelfand}. Progress in Math. 3. Birkh\"auser (1980).

\bibitem{ok-sp}Okonek, C. Spindler, H. {\em Mathematical instanton bundles on $P^{2n+1}$}, J. Reine Angew. Math., 364
(1986), pp. 3550.

\bibitem{ot-tr}Ottaviani, G. Trautmann, G. {\em  The tangent space at a special symplectic instanton bundle on $\P^{2n+1}$}. Manuscripta Math. 85, no. 1 (1994), 97-107.
 
\bibitem{qin} Qing, L. {\em Algebraic Geometry and arithmetic curves}.
 Oxford University Press, New York (2002).

\bibitem{sa}  Salamon, S.M. {\em Quaternionic structures and twistor spaces}. In T.J. Willmore
and N.J. Hitchin, editors, Global Riemannian geometry (Durham 1982), pages 6574, 1984.

\bibitem{sp-tr}	Spindler, H. Trautmann, G. {\em Special instanton bundles on $\P^{2n+1}$, their geometry and their moduli}. Math. Ann. 1990, Volume 286, Iss. 1-3, pp. 559-592.




\end{thebibliography}
\end{document}